# GLOBAL EXPONENTIAL STABILIZATION OF FREEWAY MODELS

**Iasson Karafyllis***, **Maria Kontorinaki**** and **Markos Papageorgiou****
*Dept. of Mathematics, National Technical University of Athens,
Zografou Campus, 15780, Athens, Greece (email: iasonkar@central.ntua.gr )

**Dynamic Systems and Simulation Laboratory,
Technical University of Crete, Chania, 73100, Greece
(emails: markos@dssl.tuc.gr , mkontorinaki@dssl.tuc.gr )

**Abstract**
This work is devoted to the construction of feedback laws which guarantee the robust global exponential stability of the uncongested equilibrium point for general discrete-time freeway models. The feedback construction is based on a control Lyapunov function approach and exploits certain important properties of freeway models. The developed feedback laws are tested in simulation and a detailed comparison is made with existing feedback laws in the literature. The robustness properties of the corresponding closed-loop system with respect to measurement errors are also studied.

## 1. Introduction

Freeway traffic congestion during peak periods and incidents has become a significant problem for modern societies, which leads to excessive delays, reduced traffic safety, increased fuel consumption and environmental pollution. The main traffic control measures employed to tackle traffic congestion, are ramp metering (RM) and variable speed limits (VSL). RM is implemented by use of traffic lights positioned at on-ramps to control the entering traffic flow [26]. VSL are used for speed harmonization, but recent studies have demonstrated that it may be used as a mainstream metering device as well [3]. To achieve their goal, these control measures must be driven by appropriate control strategies. A branch of related research has considered nonlinear optimal control and Model Predictive Control as a network-wide freeway traffic control approach, see, e.g. [1, 2, 9, 12]. However, possibly due to the involved control strategy complexity, none of the proposed methods has advanced to a field-operational tool. Another significant branch of freeway traffic control research has considered explicit feedback control approaches to tackle congestion problems. A pioneering development in this direction was the I-type local feedback ramp metering regulator ALINEA [24], which has been used in hundreds of successful field implementations around the world, see, e.g. [25, 27]. ALINEA controls the traffic entering from an on-ramp and targets a critical density in the mainstream merging segment so as to maximize the freeway throughput. Other proposed local feedback control algorithms for ramp metering include [13, 15, 29, 30], to mention just a few. Various extensions and modifications of ALINEA were proposed and field-implemented over the years to address specific emerging needs. Most relevant in the present context is the extension to a PI-type regulator so as to efficiently address bottlenecks which are located far downstream of the merge area [31]; and the parallel deployment of PI-type regulators to address multiple potential bottlenecks downstream of the metered on-ramp [32]. On the other hand, feedback control approaches for mainstream traffic control by use of VSL have been rather sparse, see [4]; see also [14] for a recent extension to the multiple bottleneck case.

To adequately address the increasing freeway traffic congestion problems, it is essential to investigate, develop and deploy the potentially most efficient methods; and recent control theory advances should be appropriately exploited to this end. In this work, we provide a rigorous methodology for the construction of explicit feedback laws that guarantee the robust global exponential stability of the uncongested equilibrium point for general nonlinear discrete-time freeway models. We focus on discrete-time freeway models which are generalized versions of the known first-order discrete Godunov approximations [8] to the kinematic-wave partial differential equation of the LWR-model (see [23, 28]) with nonlinear ([19]) or piecewise linear (Cell Transmission Model - CTM, [6, 7, 10]) outflow functions (fundamental diagram). Specifically, the constructed class of freeway models allows for: (i) consideration of generally nonlinear (including piecewise linear) fundamental diagrams; (ii) consideration of all possible cases for the relative priorities of the inflows at freeway nodes (see [7]), and even for time-varying and unknown priority rules; and (iii) modified demand functions according to [20], to account for the capacity drop phenomenon which is not reflected in the classical LWR-model and its Godunov discretization. The construction of the robust global exponential feedback stabilizer is based on the Control Lyapunov Function (CLF) approach (see [16]) as well as on certain important properties of freeway models. In summary, the contribution of the present work is threefold:

- a CLF is constructed for a class of freeway models; the formulas for the Lyapunov function are explicit and can be used in a straightforward way for various purposes;
- important properties of general nonlinear and uncertain discrete-time freeway models are proved;
- a parameterized family of global exponential feedback stabilizers for the uncongested equilibrium point of freeway models is constructed. The achieved stabilization is robust with respect to all priority rules that can be used for the inflows.

A comparison is made, by means of simulation, with existing feedback laws proposed in the literature and employed in practice. More specifically, we consider the Random Located Bottleneck (RLB) PI-type regulator which was proposed in [32] and is the most sophisticated of the very few comparable feedback regulators that have been employed in field operations [27]. The simulations, presented in Section 4 of the present work, study the performance of the corresponding closed-loop systems, as well as their robustness under the effect of measurement errors. It was found that the performance and the robustness properties achieved by the proposed feedback law were better than or comparable to the performance and the robustness properties induced by the RLB PI regulator. Ongoing and future work addresses further robustness issues in presence of modelling errors and persistent disturbances.

**Definitions and Notation.** Throughout this manuscript, we adopt the following notation and terminology:

* $\Re_+ := [0, +\infty)$. For every set $S$, $S^n = \underbrace{S \times \ldots \times S}_{n \text{ times}}$ for every positive integer $n$.

* By $C^0(A;\Omega)$, we denote the class of continuous functions on $A \subseteq \Re^n$, which take values in $\Omega \subseteq \Re^m$. By $C^k(A;\Omega)$, where $k \geq 1$ is an integer, we denote the class of functions on $A \subseteq \Re^n$ with continuous derivatives of order $k$, which take values in $\Omega \subseteq \Re^m$.

* Let $x \in \Re^n$. The transpose of $x \in \Re^n$ is denoted by $x'$. By $|x|$ we denote the Euclidean norm of $x \in \Re^n$.

Let $S \subseteq \Re^n$, $D \subseteq \Re^l$ be non-empty sets and consider the uncertain, discrete-time, dynamical system

$$x^+ = F(d,x), x \in S, d \in D,\qquad(1.1)$$

where $F: D \times S \to S$ is a mapping. Let $x^* \in S$ be an equilibrium point of (1.1), i.e., $x^* \in S$ satisfies $x^* = F(d,x^*)$ for all $d \in D$. Notice that $x \in S$ denotes the state of (1.1) while $d \in D$ denotes a vanishing perturbation, i.e., a disturbance that does not change the position of the equilibrium point of the system.

We use the following definitions throughout the paper.

**Definition 1.1:** *We say that $x^* \in S$ is Robustly Globally Exponentially Stable (RGES) for system (1.1) if there exist constants $M, \sigma > 0$ such that for every $x_0 \in S$ and for every sequence $\{d(t) \in D\}_{t=0}^{\infty}$ the solution $x(t)$ of (1.1) with initial condition $x(0) = x_0$ corresponding to input $\{d(t) \in D\}_{t=0}^{\infty}$ (i.e., the solution that satisfies $x(t+1) = F(d(t),x(t))$ and $x(0) = x_0$) satisfies the inequality $|x(t) - x^*| \leq M \exp(-\sigma t)|x_0 - x^*|$ for all $t \geq 0$.*

**Definition 1.2:** *A function $V: S \to \Re_+$ for which there exist constants $K_2 \geq K_1 > 0$, $p > 0$ and $\lambda \in [0,1)$ such that the inequalities $K_1|x-x^*|^p \leq V(x) \leq K_2|x-x^*|^p$ and $V(F(d,x)) \leq \lambda V(x)$ hold for all $(d,x) \in D \times S$, is called a Lyapunov function with exponent $p > 0$ for (1.1).*

**Remark 1.3:** If a Lyapunov function with exponent $p > 0$ exists for (1.1), then $x^* \in S$ is RGES. Indeed, if the state space were $\Re^n$ and not $S \subseteq \Re^n$ and if no disturbances were present, then we would be able to use Theorem 13.2 on pages 765-766 in [11]. However, since the uncertain dynamical system (1.1) is defined on $S \subseteq \Re^n$ with disturbances $d \in D$, we cannot use Theorem 13.2 on pages 765-766 in [11]. On the other hand, we can use the inequality $V(F(d,x)) \leq \lambda V(x)$ inductively and obtain the estimate $V(x(t)) \leq \lambda^t V(x(0))$ for every solution of (1.1) for every sequence $\{d(t) \in [0,1]^{n-1}\}_{t=0}^{\infty}$ and for every integer $t \geq 0$. The required exponential estimate of the solution is obtained by combining the previous estimate with the inequality $K_1|x-x^*|^p \leq V(x) \leq K_2|x-x^*|^p$.

## 2. Freeway Models and Main Result

### *2.I. Model Derivation*

We consider a freeway which consists of $n \geq 3$ components or cells; typical cell lengths may be 200-500 m. Each cell may have an external controllable inflow (e.g. from corresponding on-ramps), located near the cell's upstream boundary; and an external outflow (e.g. via corresponding off-ramps), located near the cell's downstream boundary (Figure 1). The number of vehicles at time $t \geq 0$ in component $i \in \{1,...,n\}$ is denoted by $x_i(t)$. The total outflow and the total inflow of vehicles of the component $i \in \{1,...,n\}$ at time $t \geq 0$ are denoted by $F_{i,out}(t) \geq 0$ and $F_{i,in}(t) \geq 0$, respectively. All flows during a time interval are measured in [veh]. Consequently, the balance of vehicles (conservation equation) for each component $i \in \{1,...,n\}$ gives:

$$x_i(t+1) = x_i(t) - F_{i,out}(t) + F_{i,in}(t), \; i=1,...,n, \; t \geq 0. \tag{2.1}$$

Each component of the network has storage capacity $a_i > 0$ ($i=1,...,n$). Our first assumption states that the external (off-ramp) flows from each cell are constant percentages of the total exit flow, i.e., there exist constants $p_i \in [0,1]$, $i=1,...,n$, such that:

$$\begin{pmatrix} \textit{flow of vehicles} \\ \textit{from component i to component i}+1 \end{pmatrix} = (1-p_i)F_{i,out}(t), \text{ for } i=1,...,n-1 \tag{2.2}$$

$$\begin{pmatrix} \textit{flow of vehicles from} \\ \textit{component i to regions out of the freeway} \end{pmatrix} = p_i F_{i,out}(t), \text{ for } i=1,...,n. \tag{2.3}$$

The constants $p_i$ are known as exit rates, i.e. portions of $F_{i,out}(t)$ that are bound for the off-ramp of the $i$-th cell. Since the $n$-th cell is the last downstream cell of the considered freeway, we may assume that $p_n = 1$. We also assume that $p_i < 1$ for $i=1,...,n-1$, and that all exits to regions out of the network can accommodate the respective exit flows.

Our second assumption is dealing with the attempted outflows $f_i(x_i)$, i.e. the flows that will exit the cell if there is sufficient space in the downstream cell. We assume that there exist functions $f_i \in C^0([0,a_i]; \Re_+)$ with $0 < f_i(x_i) < x_i$ for all $x_i \in (0,a_i]$ and variables $s_i(t) \in [0,1]$, $i=2,...,n$, so that:

$$F_{i-1,out}(t) = s_i(t)f_{i-1}(x_{i-1}(t)), \; i=2,...,n, \; t \geq 0 \text{ and } F_{n,out}(t) = f_n(x_n(t)). \tag{2.4}$$

The variable $s_i(t) \in [0,1]$, for each $i=2,...,n$, indicates the percentage of the attempted outflow from cell $(i-1)$ that becomes actual outflow from the same cell. The function $f_i : [0,a_i] \to \Re_+$ is called, in the specialized literature of Traffic Engineering (see, e.g., [19]), the demand-part of the fundamental diagram of the $i$-th cell, i.e. the flow that will exit the cell $i$ if there is sufficient space in the downstream cell $i+1$. Notice that equation (2.4) for $F_{n,out}(t)$ follows from our assumption that all exits to regions out of the network can accommodate the exit flows.

Let $u_i > 0$ ($i = 1,...,n$) denote the attempted external inflow to component $i \in \{1,...,n\}$ from the region out of the freeway. Typically, $u_i$, $i = 2,...,n$, correspond to external on-ramp flows which may be determined by a ramp metering control strategy. For the very first cell 1, we assume, for convenience, that there is just one external inflow, $u_1$. Let the variables $w_i(t) \in [0,1]$, $i = 1,...,n$, indicate the percentage of the attempted external inflow to component $i \in \{1,...,n\}$ that becomes actual inflow. Then, we obtain from (2.2) and (2.4):

$$F_{1,in}(t) = w_1(t)u_1(t) \text{ and } F_{i,in}(t) = w_i(t)u_i(t) + s_i(t)(1-p_{i-1})f_{i-1}(x_{i-1}(t)), \ i = 2,...,n. \quad (2.5)$$

Our next assumption is derived from the Godunov discretization (see [19]), which requires that the inflow of vehicles at the cell $i \in \{1,...,n\}$ at time $t \geq 0$, denoted by $F_{i,in}(t) \geq 0$, cannot exceed the supply function of cell $i \in \{1,...,n\}$ at time $t \geq 0$, i.e.,

$$F_{i,in}(t) \leq \min(q_i, c_i(a_i - x_i(t))), \ i = 1,...,n, \ t \geq 0 \quad (2.6)$$

where $q_i \in (0, +\infty)$ denotes the maximum flow that the $i$-th cell can receive (or the capacity flow of the $i$-th cell) and $c_i \in (0,1]$ ($i = 1,...,n$) is the jam velocity of the $i$-th cell.

Following [7], we assume that, when the total demand flow of a cell is lower than the supply of the downstream cell, i.e. when $u_i(t) + (1-p_{i-1})f_{i-1}(x_{i-1}(t)) \leq \min(q_i, c_i(a_i - x_i(t)))$ for some $i \in \{2,...,n\}$, then the demand flow can be fully accommodated by the downstream cell, and hence we have $s_i(t) = w_i(t) = 1$. Similarly, when $u_1(t) \leq \min(q_1, c_1(a_1 - x_1(t)))$, then we have $w_1(t) = 1$. In contrast, when the total demand flow of a cell is higher than the supply of the downstream cell, i.e. when $u_i(t) + (1-p_{i-1})f_{i-1}(x_{i-1}(t)) > \min(q_i, c_i(a_i - x_i(t)))$ for some $i \in \{2,...,n\}$ (or when $u_1(t) > \min(q_1, c_1(a_1 - x_1(t)))$), then the demand flow cannot be fully accommodated by the downstream cell, and the actual flow is determined by the supply function, i.e. we have $F_{i,in}(t) = \min(q_i, c_i(a_i - x_i(t)))$ (or $F_{1,in}(t) = \min(q_1, c_1(a_1 - x_1(t)))$). Therefore, we get:

$$F_{1,in}(t) = \min(q_1, c_1(a_1 - x_1(t)), u_1(t)), \ t \geq 0 \quad (2.7)$$

$$s_i(t) = (1-d_i(t))\min\left(1, \max\left(0, \frac{\min(q_i, c_i(a_i - x_i(t))) - u_i(t)}{(1-p_{i-1})f_{i-1}(x_{i-1}(t))}\right)\right) + d_i(t)\min\left(1, \frac{\min(q_i, c_i(a_i - x_i(t)))}{(1-p_{i-1})f_{i-1}(x_{i-1}(t))}\right),$$
$$i = 2,...,n, \ t \geq 0 \quad (2.8)$$

$$F_{i,in}(t) = \min(q_i, c_i(a_i - x_i(t)), u_i(t) + (1-p_{i-1})f_{i-1}(x_{i-1}(t))), \ i = 2,...,n, \ t \geq 0 \quad (2.9)$$

where

$$d_i(t) \in [0,1], \ i = 2,...,n, \ t \geq 0 \quad (2.10)$$

are time-varying parameters. Note that, if the supply is higher than the total demand, then (2.8) yields $s_i = 1$, irrespective of the value of $d_i$, since the total demand flow can be accommodated by the downstream cell. Thus, the parameter $d_i$ determines the relative inflow priorities, when the downstream supply prevails. Specifically, when $d_i(t) = 0$, then the on-ramp inflow has absolute priority over the internal inflow; on the other hand, when $d_i(t) = 1$, then the internal inflow has absolute priority over the on-ramp inflow; while intermediate values of $d_i$ reflect intermediate priority cases. The parameters $d_i(t) \in [0,1]$ are treated as unknown parameters (disturbances). Notice that by introducing the parameters $d_i(t) \in [0,1]$ (and by allowing them to be time-varying), we have taken into account all possible cases for the relative priorities of the inflows (and we also allow the priority rules to be time-varying); see [5, 17, 18, 21, 22] for freeway models with specific priority rules, which are special cases of our general approach.

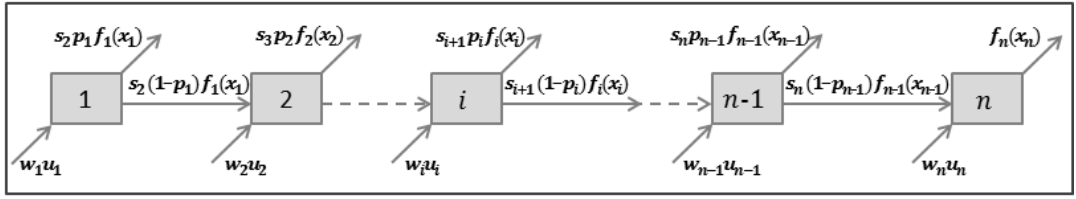

**Figure 1:** The freeway model (schematically).

All the above are illustrated in Figure 1. Combining equations (2.1), (2.2), (2.3), (2.4), (2.7) and (2.9) we obtain the following discrete-time dynamical system:

$$\begin{aligned}x_1^+ &= x_1 - s_2 f_1(x_1) + \min\left(q_1, c_1(a_1 - x_1), u_1\right) \\ &= x_1 - s_2 f_1(x_1) + w_1 u_1\end{aligned} \quad (2.11)$$

$$\begin{aligned}x_i^+ &= x_i - s_{i+1} f_i(x_i) + \min\left(q_i, c_i(a_i - x_i), u_i + (1 - p_{i-1})f_{i-1}(x_{i-1})\right) \\ &= x_i - s_{i+1} f_i(x_i) + w_i u_i + s_i(1 - p_{i-1})f_{i-1}(x_{i-1})\end{aligned}, \text{ for } i = 2,\ldots,n-1 \quad (2.12)$$

$$\begin{aligned}x_n^+ &= x_n - f_n(x_n) + \min\left(q_n, c_n(a_n - x_n), u_n + (1 - p_{n-1})f_{n-1}(x_{n-1})\right) \\ &= x_n - f_n(x_n) + w_n u_n + s_n(1 - p_{n-1})f_{n-1}(x_{n-1})\end{aligned} \quad (2.13)$$

where $s_i \in [0,1]$, $i = 2,\ldots,n$ are given by (2.8). The values of $w_i \in [0,1]$, $i = 1,\ldots,n$, may also be similarly derived from (2.5) when $u_i > 0$ but they are not needed in what follows. Define $S = (0, a_1] \times (0, a_2] \times \ldots \times (0, a_n]$. Since the functions $f_i : [0, a_i] \to \Re_+$ satisfy $0 < f_i(x_i) < x_i$ for all $x_i \in (0, a_i]$, it follows that (2.11), (2.12), (2.13) is an uncertain control system on $S$ (i.e., $x = (x_1,\ldots,x_n)' \in S$) with inputs $u = (u_1,\ldots,u_n)' \in (0,+\infty) \times \Re_+^{n-1}$ and disturbances $d = (d_2,\ldots,d_n) \in [0,1]^{n-1}$. We emphasize again that the uncertainty $d \in [0,1]^{n-1}$ appears in the equations (2.11), (2.12) and (2.13) only when the supply function prevails, i.e., only when $u_i(t) + (1 - p_{i-1})f_{i-1}(x_{i-1}(t)) > \min\left(q_i, c_i(a_i - x_i(t))\right)$ for some $i \in \{2,\ldots,n\}$.

We make the following assumption for the functions $f_i : [0, a_i] \to \Re_+$ ($i = 1,\ldots,n$):

**(H)** *The function $f_i \in C^0([0, a_i]; \Re_+)$ satisfies $0 < f_i(z) < z$ for all $z \in (0, a_i]$. There exists $\delta_i \in (0, a_i]$ such that $f_i$ is increasing on $[0, \delta_i]$ and non-increasing on $[\delta_i, a_i]$. Moreover, there exist constants $L_i \in (0,1)$, $\tilde{\delta}_i \in (0, \delta_i]$ such that $f_i : [0, a_i] \to \Re_+$ is $C^1$ on $(0, \delta_i)$, $1 - L_i \leq f_i'(z)$ for all $z \in (0, \tilde{\delta}_i)$, $f_i'(z) \leq 1$ for all $z \in (0, \delta_i)$.*

Assumption (H) reflects the basic properties of the so-called "demand function" [19] in the Godunov discretization; whereby $\delta_i$ is the critical density, where $f_i(x_i)$ achieves a maximum value. In other words, the fundamental diagram (FD) of cell $i$ is composed by the increasing function $f_i(x_i)$ for $x_i \in [0, \delta_i]$; and by the non-increasing supply function $\min(q_i, c_i(a_i - x_i))$ for $x_i \in [\delta_i, a_i]$. Note, however, that Assumption (H) includes the possibility of reduced demand flow for overcritical densities (i.e., when $x_i(t) \geq \delta_i$), since $f_i(x_i)$ is allowed to be decreasing for $x_i \in [\delta_i, a_i]$; this could be used to reflect the capacity drop phenomenon as proposed in [20]. In conclusion, the model (2.11)-(2.13) is a generalized version of the known first-order discrete Godunov approximation to the kinematic-wave partial differential equation of the LWR-model (see [23, 28]) with nonlinear ([19]) or piecewise linear (Cell Transmission Model - CTM, [6, 7]) fundamental diagram. However, the presented framework can also accommodate recent modifications of the LWR-model as in [20] to reflect the so-called capacity drop phenomenon. Notice that the piecewise smooth selections $f_i(z) = q_i \delta_i^{-1} z$ for $z \in [0, \delta_i]$ and $f_i(z) = q_i$ for $z \in (\delta_i, a_i]$ ($i = 1, ..., n$) with $a_i \geq \delta_i + c_i^{-1} q_i$ allow us to obtain the standard CTM with: (i) triangular-FD (if $a_i = \delta_i + c_i^{-1} q_i$); and (ii) trapezoidal-FD (if $a_i > \delta_i + c_i^{-1} q_i$). In the latter case, assumption (H) holds with arbitrary $\tilde{\delta}_i \in (0, \delta_i]$.

Define the vector field $\tilde{F} : [0,1]^{n-1} \times S \times (0, +\infty) \times \mathfrak{R}_+^{n-1} \to S$ for all $x \in S := (0, a_1] \times ... (0, a_n]$, $d = (d_2, ..., d_n) \in D = [0,1]^{n-1}$ and $u = (u_1, ..., u_n) \in (0, +\infty) \times \mathfrak{R}_+^{n-1}$:

$$\tilde{F}(d, x, u) = (\tilde{F}_1(d, x, u), ..., \tilde{F}_n(d, x, u))' \in \mathfrak{R}^n$$
$$\text{with } \tilde{F}_1(d, x, u) := x_1 - s_2 f_1(x_1) + \min(q_1, c_1(a_1 - x_1), u_1),$$
$$\tilde{F}_i(d, x, u) = x_i - s_{i+1} f_i(x_i) + \min(q_i, c_i(a_i - x_i), u_i + (1 - p_{i-1}) f_{i-1}(x_{i-1})), \text{ for } i = 2, ..., n-1,$$
$$\tilde{F}_n(d, x, u) = x_n - f_n(x_n) + \min(q_n, c_n(a_n - x_n), u_n + (1 - p_{n-1}) f_{n-1}(x_{n-1})) \text{ and}$$
$$s_i = (1 - d_i) \min\left(1, \max\left(0, \frac{\min(q_i, c_i(a_i - x_i)) - u_i}{(1 - p_{i-1}) f_{i-1}(x_{i-1})}\right)\right) + d_i \min\left(1, \frac{\min(q_i, c_i(a_i - x_i))}{(1 - p_{i-1}) f_{i-1}(x_{i-1})}\right), \text{ for } i = 2, ..., n.$$

(2.14)

Notice that, using definition (2.14), the control system (2.11), (2.12), (2.13) can be written in the following vector form:

$$x^+ = \tilde{F}(d, x, u)$$
$$x \in S, d \in D, u \in (0, +\infty) \times \mathfrak{R}_+^{n-1} \quad (2.15)$$

## 2.II. Main Result

Consider the freeway model (2.15) under assumption (H). We suppose that there exist $u_1^* > 0$, $u_i^* \geq 0$ ($i = 2, ..., n$) and a vector $x^* = (x_1^*, ..., x_n^*) \in (0, \tilde{\delta}_1) \times ... (0, \tilde{\delta}_n)$ with

$$f_1(x_1^*) = u_1^*, \quad f_i(x_i^*) = u_i^* + (1 - p_{i-1}) f_{i-1}(x_{i-1}^*) = u_i^* + \sum_{j=1}^{i-1} \left( \prod_{k=j}^{i-1} (1 - p_k) \right) u_j^* \quad (i = 2, ..., n) \quad (2.16)$$

and $u_1^* < \min\left(q_1, c_1(a_1 - x_1^*)\right)$, $u_i^* + (1 - p_{i-1}) f_{i-1}(x_{i-1}^*) < \min\left(q_i, c_i(a_i - x_i^*)\right)$ ($i = 2, ..., n$). (2.17)

This is the uncongested equilibrium point of the freeway model (2.15). Notice that assumption (H) guarantees that an uncongested equilibrium point always exists for the freeway model (2.15) when $u_1^* > 0$ and $u_i^* \geq 0$ ($i = 2,...,n$) are sufficiently small. The uncongested equilibrium point is not globally exponentially stable for arbitrary $u_1^* > 0$, $u_i^* \geq 0$ ($i = 2,...,n$); indeed, for relatively large values of external demands $u_1^* > 0$, $u_i^* \geq 0$ ($i = 2,...,n$) there exist other equilibria for model (2.15) (congested equilibria) for which the cell densities are large and can attract the solution of (2.15) (see the numerical example 4.2 in the section 4).

The following result is our main result in feedback design. The result shows that a continuous, robust, global exponential stabilizer exists for every freeway model of the form (2.15) under assumption (H).

**Theorem 2.1:** *Consider system (2.15) with $n \geq 3$ under assumption (H). Then there exist a subset $R \subseteq \{1,...,n\}$ of the set of all indices $i \in \{1,...,n\}$ with $u_i^* > 0$, constants $\sigma \in (0,1]$, $b_i \in (0, u_i^*)$ for $i \in R$ and a constant $\tau^* > 0$ such that for every $\tau \in (0, \tau^*)$ the feedback law $k : S \to \Re_+^n$ defined by:*

$$k(x) := (k_1(x),...,k_n(x))' \in \Re^n \text{ with}$$

$$k_i(x) := \max\left(u_i^* - \gamma_i \Xi(x), b_i\right), \text{ for all } x \in S, i \in R \text{ and } k_i(x) := u_i^*, \text{ for all } x \in S, i \notin R \quad (2.18)$$

*where $\gamma_i := \tau^{-1}(u_i^* - b_i)$ and*

$$\Xi(x) := \sum_{i=1}^{n} \sigma^i \max\left(0, x_i - x_i^*\right), \text{ for all } x \in S \quad (2.19)$$

*achieves robust global exponential stabilization of the uncongested equilibrium point $x^*$ of system (2.15), i.e., $x^*$ is RGES for the closed-loop system (2.15) with $u = k(x)$. Moreover, for every $\tau \in (0, \tau^*)$, there exist constants $Q, h, \theta, A, K > 0$ so that the function $V : S \to \Re_+$ defined by:*

$$V(x) := \sum_{i=1}^{n} \sigma^i \left|x_i - x_i^*\right| + A\Xi(x) + K \max\left(0, \sum_{i=1}^{n} I_i(x) - P(x)\right), \text{ for all } x \in S \quad (2.20)$$

*where $I_j(x) := \sum_{i=1}^{j} x_i$ for $j = 1,...,n$ and*

$$P(x) := Q - \theta \min\left(h, \Xi(x)\right) \quad (2.21)$$

*is a Lyapunov function with exponent 1 for the closed-loop system (2.15) with $u = k(x)$.*

Although Theorem 2.1 is an existence result, its proof is constructive and provides formulae for all constants and for the index set $R$ (see following sections). Notice that the index set $R$ is the set of all inflows that must be controlled in order to be able to guarantee that the uncongested equilibrium point is RGES; consequently, the knowledge of the index set $R$ is critical.

The importance of Theorem 2.1 lies on the facts that:

- Theorem 2.1 provides a family of robust global exponential stabilizers (parameterized by the parameter $\tau \in (0, \tau^*)$) and an explicit formula for the feedback law (formula (2.18));
- the achieved stabilization result is robust for all possible (and even time-varying) priority rules for the junctions that may apply at specific freeways; thus, there is no need to know or estimate the applied priority rules;
- Theorem 2.1 provides an explicit formula for the Lyapunov function of the closed-loop system. This is important, because the knowledge of the Lyapunov function can allow the study of the robustness of the closed-loop system to various disturbances (measurement errors, modeling errors, etc.) as well for the study of the effect of interconnections of freeways (by means of the small-gain theorem; see [16]).

The main idea behind the proof of Theorem 2.1 is the construction of the Lyapunov function of the closed-loop system, which acts as a Control Lyapunov Function (CLF; see [16]) for the open-loop system. The construction of the Lyapunov function is based on the observation that there are no congestion phenomena when the cell densities are sufficiently small, i.e.,

"There exists a set $\Omega \subset S$ of the form $\Omega = (0, \mu_1] \times ... \times (0, \mu_n]$, where $\mu_i > 0$ for $i = 1,...,n$ are constants, such that no congestion phenomena are present when $x \in \Omega$."

The existence of the set $\Omega \subset S$ is important because, when no congestion phenomena are present, then the freeway model admits the simple (cascade) form:

$$x_1^+ = x_1 - f_1(x_1) + u_1, \quad x_i^+ = x_i - f_i(x_i) + (1 - p_{i-1})f_{i-1}(x_{i-1}) + u_{i+1}, \text{ for } i = 2,...,n, \; x \in \Omega$$

and a Lyapunov function for the above form can be a function of the form $V_1(x) := \sum_{i=1}^{n} \sigma^i |x_i - x_i^*| + A\Xi(x)$, where $\Xi(x) := \sum_{i=1}^{n} \sigma^i \max(0, x_i - x_i^*)$; and $\sigma \in (0,1]$ and $A > 0$ are appropriate constants. The Lyapunov function for the freeway model is the linear combination of the "Lyapunov function" for the uncongested model (i.e., $V_1(x) := \sum_{i=1}^{n} \sigma^i |x_i - x_i^*| + A\Xi(x)$) and a penalty term, i.e., the term $\max\left(0, \sum_{i=1}^{n} I_i(x) - P(x)\right)$, that penalizes large cell densities (and thus penalizes the possibility of the state being out of the set $\Omega \subset S$). The appropriate selection of the weight of the penalty term $K > 0$ forces the selected control action to lead the state in the set $\Omega \subset S$ (see Figure 2). In other words, the construction of the CLF guarantees that the control action will first eliminate all congestion phenomena and then will drive the state to the desired equilibrium.

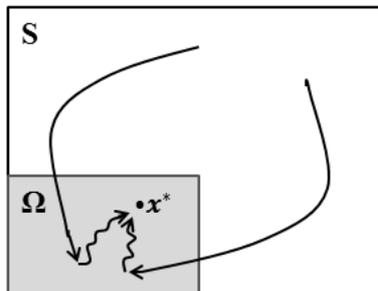

**Figure 2:** Idea behind Theorem 2.1

## 3. Proof of Main Result

Assumption (H) has non-trivial consequences. A list of the most important consequences of assumption (H) is given below. All following consequences are exploited in the proof of Theorem 2.1.

Consequences of assumption (H):

**(C1)** *The mappings $[0, a_i] \ni z \to (z - f_i(z)) \geq 0$ are non-decreasing for $i = 1,...,n$.*

Property (C1) is a direct consequence of the fact that $f_i'(z) \leq 1$ for all $z \in (0, \delta_i)$ and the fact that $f_i$ is non-increasing on $[\delta_i, a_i]$.

**(C2)** *For each $i = 1,...,n$ there exist constants $\lambda_i \in (0,1)$, $G_i \in [0,1]$ such that:*

$$\left| x_i - x_i^* - f_i(x_i) + f_i(x_i^*) \right| \leq \lambda_i \left| x_i - x_i^* \right| \text{ and } \left| f_i(x_i) - f_i(x_i^*) \right| \leq G_i \left| x_i - x_i^* \right|, \text{ for all } x_i, x_i^* \in [0, \tilde{\delta}_i]. \quad (3.1)$$

Property (C2) is a direct consequence of the fact that there exist constants $L_i \in (0,1)$, $\tilde{\delta}_i \in (0, \delta_i]$ such that $f_i : [0, a_i] \to \Re_+$ is $C^1$ on $(0, \delta_i)$ and $1 - L_i \leq f_i'(z)$ for all $z \in (0, \tilde{\delta}_i)$ and $f_i'(z) \leq 1$ for all $z \in (0, \delta_i)$. We conclude that (3.1) holds with $\lambda_i = L_i \in (0,1)$ and $G_i = 1$.

**(C3)** *There exist constants $\theta_i > 0$ ($i = 1,...,n$) such that $f_i(z) \geq \theta_i z$ for all $z \in [0, a_i]$ and $i = 1,...,n$.*

Property (C3) is a direct consequence of the fact that $f_i(z) \geq (1 - L_i)z$ for all $z \in (0, \tilde{\delta}_i]$ (a direct consequence of the fact that $1 - L_i \leq f_i'(z)$ for all $z \in (0, \tilde{\delta}_i)$) and the fact that $0 < f_i(z) < z$ for all $z \in (0, a_i]$ (for example, the selection $\theta_i = \min\left(1 - L_i, \min_{z \in [\tilde{\delta}_i, a_i]} \left(z^{-1} f_i(z)\right)\right) > 0$ satisfies (C3)).

The following consequence provides a useful linear lower bound on a weighted sum of exit rates. Its proof is provided in the Appendix.

**(C4)** *For every $r_2,...,r_n \geq 0$ with $r_i < \min(q_i, c_i a_i)$ for $i = 2,...,n$, there exists a constant $C > 0$ such that the following inequality holds for all $x \in S := (0, a_1] \times ... (0, a_n]$, $u \in U = (0, +\infty) \times [0, r_2] \times ... \times [0, r_n]$, $d = (d_2,...,d_n) \in [0,1]^{n-1}$ with $p_n = 1 = s_{n+1}$:*

$$\sum_{i=1}^{n} (1 + p_i(n-i)) s_{i+1} f_i(x_i) \geq C \sum_{i=1}^{n} (n + 1 - i) x_i. \quad (3.2)$$

**Remark 3.1:** The proof of property (C4) implies that the constant $C > 0$ can be estimated in a straightforward way. We define the positive constants $Y_i > 0$ ($i = 1,...,n$) using the recursive formula:

$$Y_{k-1} = \min\left(Y_k, \frac{1+p_{k-1}(n+1-k)}{n+2-k} l_k \theta_{k-1}, \frac{(n+1-k)\left(a_k - c_k^{-1} r_k\right) Y_k}{2(n+1-k)a_k + 2(n+2-k)a_{k-1}}, \frac{q_k - r_k}{1-p_{k-1}} \frac{1+p_{k-1}(n+1-k)}{(n+2-k)a_{k-1}}\right) \quad (3.3)$$

for $k = n, n-1, ..., 2$, with $Y_n = \theta_n$, where $\theta_i > 0$ ($i = 1, ..., n$) are the constants involved in Property (C3) and $l_k = \min\left(1, \frac{c_k a_k - r_k}{2(1-p_{k-1})f_{k-1}(\delta_{k-1})}\right)$ for $k = 2, ..., n$. Then the constant $C > 0$ can be selected as $C = Y_1$. However, the estimation of the constant $C > 0$ by the recursive formula (3.3) is, in general, conservative.

Finally, the last consequence provides useful equalities and inequalities for a weighted sum of all vehicle densities of the freeway. Its proof is provided in the Appendix.

**(C5)** *The following equality holds for all* $x \in S := (0, a_1] \times ... (0, a_n]$, $u = (u_1, ..., u_n)' \in (0, +\infty) \times \mathfrak{R}_+^{n-1}$, $d = (d_2, ..., d_n) \in [0,1]^{n-1}$ *with* $p_n = 1 = s_{n+1}$:

$$\sum_{i=1}^n I_i(x^+) = \sum_{i=1}^n I_i(x) + \sum_{i=1}^n (n+1-i) w_i u_i - \sum_{i=1}^n (1 + p_i(n-i)) s_{i+1} f_i(x_i) \quad (3.4)$$

*where* $I_j(x) := \sum_{i=1}^j x_i$ *for* $j = 1, ..., n$. *Moreover, for every* $r_2, ..., r_n \geq 0$ *with* $r_i < \min(q_i, c_i a_i)$ *for* $i = 2, ..., n$, *the following inequality holds:*

$$\sum_{i=1}^n I_i(x^+) \leq (1-C)\sum_{i=1}^n I_i(x) + \sum_{i=1}^n (n+1-i) u_i, \text{ for all } (x, u, d) \in S \times U \times [0,1]^{n-1} \quad (3.5)$$

*where* $U = (0, +\infty) \times [0, r_2] \times ... \times [0, r_n]$, *and* $C > 0$ *is the constant involved in (3.2).*

We are now ready to provide the proof of Theorem 2.1.

**Proof of Theorem 2.1:** Define $\beta_n = \tilde{\delta}_n$ and $\beta_i \in (0, \tilde{\delta}_i]$, for $i = 1, ..., n-1$ to be the unique solution of the equation

$$f_i(\beta_i) = \min\left(f_i(\tilde{\delta}_i), \frac{q_{i+1} - u_{i+1}^*}{1 - p_i}\right). \quad (3.6)$$

Due to the inequalities (2.17) and the fact that $\beta_n = \tilde{\delta}_n$, it follows that $\beta_i > x_i^*$ for $i = 1, ..., n$. Define $\omega_i = c_i\left(a_i - x_i^*\right) - u_i^* - (1 - p_{i-1})f_{i-1}(x_{i-1}^*)$ for $i = 2, ..., n$ and $\omega_1 = c_1\left(a_1 - x_1^*\right) - u_1^*$. Next define:

$$\mu_i = \min\left(\beta_i, x_i^* + \frac{\omega_i}{2c_i}, x_i^* + \frac{\omega_{i+1}}{2(1-p_i)}\right) \text{ for } i = 1, ..., n-1 \text{ and } \mu_n = \min\left(\beta_n, x_n^* + \frac{\omega_n}{2c_n}\right) \quad (3.7)$$

Again, due to the inequalities (2.17) and the fact that $\beta_i > x_i^*$ for $i = 1, ..., n$, it follows that $\mu_i > x_i^*$ for $i = 1, ..., n$.

It follows from (2.5), (2.7), (2.8), (2.9), (2.11), (2.12), (2.13), (3.6) and (3.7) that the following equations hold when $x \in \Omega = (0, \mu_1] \times ... \times (0, \mu_n]$ and $u_i \in [0, u_i^*]$ for $i = 1,...,n$:

$$w_i = 1, \text{ for } i = 1,...,n \text{ and } s_i = 1, \text{ for } i = 2,...,n \tag{3.8}$$

$$x_1^+ = x_1 - f_1(x_1) + u_1, \; x_i^+ = x_i - f_i(x_i) + (1 - p_{i-1})f_{i-1}(x_{i-1}) + u_{i+1}, \text{ for } i = 2,...,n. \tag{3.9}$$

To see this, notice that for all $x \in \Omega = (0, \mu_1] \times ... \times (0, \mu_n]$ and $u \in [0, u_1^*] \times ... \times [0, u_n^*]$ we have $u_1 \leq \min(q_1, c_1(a_1 - x_1))$ and $u_i + (1 - p_{i-1})f_{i-1}(x_{i-1}) \leq \min(q_i, c_i(a_i - x_i))$ for $i = 2,...,n$. Indeed, assumption (H) in conjunction with equation (3.6) and definition (3.7) implies that

$$u_i + (1 - p_{i-1})f_{i-1}(x_{i-1}) \leq u_i^* + (1 - p_{i-1})f_{i-1}(\mu_{i-1}) \leq u_i^* + (1 - p_{i-1})f_{i-1}(\beta_{i-1}) \leq q_i$$

for $i = 2,...,n$ and for all $x \in \Omega = (0, \mu_1] \times ... \times (0, \mu_n]$ and $u \in [0, u_1^*] \times ... \times [0, u_n^*]$. The inequality $u_1 \leq u_1^* \leq q_1$ is directly implied by (2.17). Moreover, assumption (H) (and particularly the fact that $f_i$ is increasing on $[0, \delta_i]$ with $f_i'(z) \leq 1$ for all $z \in (0, \delta_i)$ for $i = 1,...,n$) in conjunction with (2.17) and definition (3.7) implies that

$$u_i + (1 - p_{i-1})f_{i-1}(x_{i-1}) \leq u_i^* + (1 - p_{i-1})f_{i-1}(x_{i-1}^*) + (1 - p_{i-1})\left(f_{i-1}(x_{i-1}) - f_{i-1}(x_{i-1}^*)\right)$$
$$\leq u_i^* + (1 - p_{i-1})f_{i-1}(x_{i-1}^*) + (1 - p_{i-1})\max\left(0, x_{i-1} - x_{i-1}^*\right)$$
$$\leq u_i^* + (1 - p_{i-1})f_{i-1}(x_{i-1}^*) + \frac{\omega_i}{2} = \frac{1}{2}c_i\left(a_i - x_i^*\right) + \frac{1}{2}u_i^* + \frac{1}{2}(1 - p_{i-1})f_{i-1}(x_{i-1}^*) = c_i\left(a_i - x_i^*\right) - \frac{\omega_i}{2}$$
$$\leq c_i\left(a_i - x_i^*\right) - c_i \max\left(0, x_i - x_i^*\right) \leq c_i\left(a_i - x_i^*\right) + c_i\left(x_i^* - x_i\right) = c_i\left(a_i - x_i\right)$$

for $i = 2,...,n$ and for all $x \in \Omega = (0, \mu_1] \times ... \times (0, \mu_n]$ and $u \in [0, u_1^*] \times ... \times [0, u_n^*]$. The inequality $u_1 \leq u_1^* \leq c_1(a_1 - x_1)$ is a consequence of (2.17), definition (3.7) and the inequalities

$$u_1^* \leq \frac{1}{2}u_1^* + \frac{1}{2}c_1\left(a_1 - x_1^*\right) = c_1\left(a_1 - x_1^*\right) - \frac{\omega_1}{2} \leq c_1\left(a_1 - x_1^*\right) - c_1 \max\left(0, x_1 - x_1^*\right) \leq$$
$$c_1\left(a_1 - x_1^*\right) + c_1\left(x_1^* - x_1\right) = c_1\left(a_1 - x_1\right)$$

Let $\lambda_i \in (0,1)$, $G_i \in [0,1]$ ($i = 1,...,n$), be the constants involved in Property (C2). Let $\sigma \in (0,1]$ be a constant so that

$$L := \max\left(\lambda_n, \max_{i=1,...,n-1}(\lambda_i + \sigma G_i(1 - p_i))\right) < 1. \tag{3.10}$$

Notice that $\max\left(\lambda_n, \max_{i=1,...,n-1}(\lambda_i + \sigma G_i(1 - p_i))\right) < 1$ for all $\sigma \in (0,1)$. In what follows, we have $p_n = 1 = s_{n+1}$. Let $r_i = u_i^*$ for $i = 2,...,n$ and let $C > 0$ be the constant involved in (3.2). Let $R \subseteq \{1,...,n\}$ be a subset of the set of all indices $i \in \{1,...,n\}$ for which $u_i^* > 0$ and such that:

$$\sum_{i \notin R}(n+1-i)u_i^* < \min_{i=1,\ldots,n}\left(((n-i)p_i+1)f_i(x_i^*)\right) \text{ and } \sum_{i \notin R}(n+1-i)u_i^* < C\min_{i=1,\ldots,n}\left((n+1-i)\mu_i\right) \quad (3.11)$$

where $\mu_i > x_i^*$ for $i=1,\ldots,n$ are the constants defined by (3.6). Such a set $R \subseteq \{1,\ldots,n\}$ always exists (for example, $R \subseteq \{1,\ldots,n\}$ can be the set of all indices $i \in \{1,\ldots,n\}$ for which $u_i^* > 0$). Inequalities (3.11) imply that there exist constants $\varepsilon \in (0,1)$ and $b_i \in (0, u_i^*)$ for $i \in R$ such that:

$$\sum_{i \in R}(n+1-i)b_i + \sum_{i \notin R}(n+1-i)u_i^* \leq \min_{i=1,\ldots,n}\left(((n-i)p_i+1)f_i(x_i^*)\right) \text{ and}$$
$$\sum_{i \in R}(n+1-i)b_i + \sum_{i \notin R}(n+1-i)u_i^* \leq \varepsilon C \min_{i=1,\ldots,n}\left((n+1-i)\mu_i\right). \quad (3.12)$$

We next define the following parameters:

- Define $h := \min_{i=1,\ldots,n}\left(\sigma^i(\mu_i - x_i^*)\right)$.
- Define
$$Q := \max\left(\min_{i=1,\ldots,n}\left(\mu_i(n+1-i)\right), (1-C)\sum_{i=1}^n I_i(x^*) + (1-C)h\max_{i=1,\ldots,n}\left((n+1-i)\sigma^{-i}\right) + \sum_{i=1}^n(n+1-i)u_i^*\right)$$
and $\theta := h^{-1}\left(Q - \varepsilon \min_{i=1,\ldots,n}\left((n+1-i)\mu_i\right)\right)$.
- Define $\tau^* := \min\left(h, (\theta L)^{-1}\sum_{i \in R}(n+1-i)(u_i^* - b_i)\right)$ and let $\tau \in (0, \tau^*)$.
- Define $A := 1 + (1-L)^{-1}\sum_{i \in R}\sigma^i \gamma_i$, where $\gamma_i := \tau^{-1}(u_i^* - b_i)$ for $i \in R$.
- Define $K := \dfrac{\sum_{i=1}^n \sigma^i \max(a_i - x_i^*, x_i^*) + A\sum_{i=1}^n \sigma^i(a_i - x_i^*) - (A+L)h}{(1-\varepsilon)C\min_{i=1,\ldots,n}\left((n+1-i)\mu_i\right)}$.

We next prove the implication:

If $x \in \Omega = (0, \mu_1] \times \ldots \times (0, \mu_n]$, $d \in [0,1]^{n-1}$ and $u \in [0, u_1^*] \times \ldots \times [0, u_n^*]$ then $\Xi(x^+) \leq L\Xi(x)$ (3.13)

where $L \in (0,1)$ is defined by (3.10) and $x^+ = \tilde{F}(d, x, u)$. Indeed, using (3.9) and definition (2.19), we get for all $x \in \Omega = (0, \mu_1] \times \ldots \times (0, \mu_n]$, $d \in [0,1]^{n-1}$ and $u \in [0, u_1^*] \times \ldots \times [0, u_n^*]$:

$$\Xi(x^+) = \sum_{i=2}^n \sigma^i \max\left(0, x_i - f_i(x_i) + (1-p_{i-1})f_{i-1}(x_{i-1}) + u_i - x_i^*\right) + \sigma \max\left(0, x_1 - f_1(x_1) + u_1 - x_1^*\right)$$
$$\leq \sum_{i=1}^n \sigma^i \max\left(0, x_i - f_i(x_i) + f_i(x_i^*) - x_i^*\right) + \sum_{i=2}^n \sigma^i (1 - p_{i-1})\max\left(0, f_{i-1}(x_{i-1}) - f_{i-1}(x_{i-1}^*)\right)$$
(3.14)

Using (3.1), the fact that $\mu_i \leq \tilde{\delta}_i$ for $i=1,...,n$ (a consequence of (3.6) and (3.7)) and the fact that $f_i$ is increasing on $[0, \tilde{\delta}_i]$ for $i=1,...,n$ (a consequence of assumption (H)), we get:

$$\max\left(0, f_i(x_i) - f_i(x_i^*)\right) \leq G_i \max\left(0, x_i - x_i^*\right), \text{ for all } x_i \in [0, \mu_i], \ i=1,...,n. \quad (3.15)$$

Using properties (C1), (C2) and the fact that $\mu_i \leq \tilde{\delta}_i$ for $i=1,...,n$ (a consequence of (3.6) and (3.7)), we get:

$$\max\left(0, x_i - f_i(x_i) + f_i(x_i^*) - x_i^*\right) \leq \lambda_i \max\left(0, x_i - x_i^*\right), \text{ for all } x_i \in [0, \mu_i], \ i=1,...,n. \quad (3.16)$$

Combining (3.10), (3.14), (3.15), (3.16), we obtain implication (3.13).

Next, we show the implication:

$$\text{If } x \in S, \ d \in [0,1]^{n-1} \text{ and } u \in [0, u_1^*] \times ... \times [0, u_n^*] \text{ then } P(x^+) \geq P(x). \quad (3.17)$$

where $x^+ = \tilde{F}(d, x, u)$. Indeed, (3.17) is a direct consequence of (3.13) and definition (2.21) when $x \in \Omega = (0, \mu_1] \times ... \times (0, \mu_n]$. On the other hand, when $x \in S \setminus \Omega$ there exists at least one $i \in \{1,...,n\}$ for which $x_i > \mu_i$. Therefore, definition (2.19) implies $\Xi(x) > \min_{i=1,...,n}\left(\sigma^i \left(\mu_i - x_i^*\right)\right)$, and consequently definition (2.21) gives $P(x) = Q - \theta h$ (a consequence of the fact that $h = \min_{i=1,...,n}\left(\sigma^i \left(\mu_i - x_i^*\right)\right)$). Since $P(x) \geq Q - \theta h$ for all $x \in S$ (a consequence of (2.21)), we get $P(x^+) \geq Q - \theta h = P(x)$ when $x \in S \setminus \Omega$.

In what follows, we have $x^+ = \tilde{F}(d, x, k(x))$. Next we make the following claims. Their proofs can be found in the Appendix.

**(Claim 1):** *For all $x \in S$, $d \in [0,1]^{n-1}$, the following inequality holds:*

$$V(x^+) \leq V(x) - (1-L)\sum_{i=1}^{n} \sigma^i \left|x_i - x_i^*\right|. \quad (3.18)$$

**(Claim 2):** *There exist constants $K_2 \geq K_1 > 0$ such that the following inequality holds.*

$$K_1 \left|x - x^*\right| \leq V(x) \leq K_2 \left|x - x^*\right| \text{ for all } x \in S. \quad (3.19)$$

Using (3.18), the fact that $\sigma \in (0,1]$, and (3.19), we get for all $x \in S$, $d \in [0,1]^{n-1}$:

$$V(x^+) \leq V(x) - (1-L)\sum_{i=1}^{n} \sigma^i \left|x_i - x_i^*\right| \leq V(x) - (1-L)\sigma^n \left|x - x^*\right| \leq \left(1 - (1-L)\sigma^n K_2^{-1}\right) V(x).$$

The above inequality implies that the inequality

$$V(\tilde{F}(d, x, k(x))) \leq \tilde{L} V(x) \text{ for all } x \in S, \ d \in [0,1]^{n-1} \quad (3.20)$$

holds with $\tilde{L} := 1 - (1-L)\sigma^n K_2^{-1}$. Notice that $\tilde{L} \in (0,1)$. Inequalities (3.19) and (3.20) show that the function $V : S \to \Re_+$ is a Lyapunov function with exponent 1 for the closed-loop system (2.15) with $u = k(x)$. Remark 1.3 guarantees that $x^*$ is RGES for the closed-loop system (2.15) with $u = k(x)$. The proof is complete. ◁

**Remark 3.2:** The proof of Theorem 2.1 provides a methodology for obtaining an estimation of the set $R \subseteq \{1,...,n\}$, the constant $\sigma \in (0,1]$ and the critical constant $\tau^* > 0$. Let $r_i = u_i^*$ for $i = 2,...,n$ and let $C > 0$ be the constant involved in (3.2). Select $R \subseteq \{1,...,n\}$ to be a subset of the set of all indices $i \in \{1,...,n\}$, for which $u_i^* > 0$ and for which there exist $b_i \in (0, u_i^*)$ such that:

$$\sum_{i \in R}(n+1-i)b_i + \sum_{i \notin R}(n+1-i)u_i^* \leq \min_{i=1,...,n}\left(\left((n-i)p_i + 1\right)f_i(x_i^*)\right) \text{ and}$$

$$\sum_{i \in R}(n+1-i)b_i + \sum_{i \notin R}(n+1-i)u_i^* < C \min_{i=1,...,n}\left((n+1-i)\mu_i\right)$$

where $\mu_i > x_i^*$, for $i = 1,...,n$, are the constants defined by (3.6). Let $\varepsilon \in (0,1)$ be a constant which satisfies $\sum_{i \in R}(n+1-i)b_i + \sum_{i \notin R}(n+1-i)u_i^* \leq \varepsilon C \min_{i=1,...,n}\left((n+1-i)\mu_i\right)$. The estimation of the critical constant $\tau^* > 0$ may be done in the following way:

- Select $\sigma \in (0,1]$ so that $L = \max\left(\lambda_n, \max_{i=1,...,n-1}(\lambda_i + \sigma G_i(1-p_i))\right) < 1$, where $\lambda_i \in (0,1)$, $G_i \in [0,1]$ ($i = 1,...,n$), are the constants involved in Property (C2).
- Define $h := \min_{i=1,...,n}\left(\sigma^i(\mu_i - x_i^*)\right)$.
- Define
$$Q := \max\left(\min_{i=1,...,n}(\mu_i(n+1-i)), (1-C)\sum_{i=1}^n I_i(x^*) + (1-C)h \max_{i=1,...,n}\left((n+1-i)\sigma^{-i}\right) + \sum_{i=1}^n(n+1-i)u_i^*\right)$$

and $\theta := h^{-1}\left(Q - \varepsilon \min_{i=1,...,n}\left((n+1-i)\mu_i\right)\right)$, where $I_j(x) := \sum_{i=1}^j x_i$ for $j = 1,...,n$.

The estimated value of $\tau^* > 0$ is given by $\tau^* := \min\left(h, (\theta L)^{-1}\sum_{i \in R}(n+1-i)\left(u_i^* - b_i\right)\right)$. However, the estimated value of $\tau^* > 0$, which is obtained by applying the above methodology, may be conservative (significantly smaller than the actual value).

## 4. Illustrative Examples

The issue of the selection of which specific inflows must be controlled for the stabilization of the uncongested equilibrium point of a freeway is crucial. The following example illustrates how Theorem 2.1 can be used for such a selection.

**Example 4.1:** Consider a freeway stretch, which consists of $n = 4$ cells, where the first and the third cell have one on-ramp, while there are no intermediate off-ramps (i.e., $p_i = 0$ for $i = 1, 2, 3$). Each cell has jam density $a_i = 80$ and is characterized by the same demand functions which are given by:

$$f_i(z) = \begin{cases} (9/10)z & \text{for } z \in [0, 40] \\ (9/10)(80-z) & \text{for } z \in (40, 370/9] \\ 35 & \text{for } z \in (370/9, 80] \end{cases} \quad (i = 1,...,4).$$

Assumptions (H) hold with $\delta_i = \tilde{\delta}_i = 40$, $c_i = 9/10$, $q_i = 36$ (leading to a triangular-shaped FD) and $L_i = 1/10$ ($i = 1,...,4$). Property (C3) holds with $\theta_i = 7/16$ for $i = 1,...,4$. The value of constant $C > 0$ that satisfies (3.2) was estimated by applying the procedure described in Remark 3.1 with $r_1 = 35.9$, $r_2 = r_4 = 0$ and $r_3 = 1$. It was found that $C \geq 0.005$.

Next, we consider inflows $u_1^* = 35.5$, $u_2^* = u_4^* = 0$ and $u_3^* > 0$. The uncongested equilibrium point exists for all $u_3^* < 0.5$. However, for constant inflow $u_1 = 35.5$, the uncongested equilibrium point is not globally exponentially stable due to the existence of congested equilibria. Therefore, there is a need for controlling the main inflow $u_1$. At this point, the following question becomes crucial:

"For what values of $u_3^* \in (0, 0.5)$ can the uncongested equilibrium point be globally exponentially stabilized by controlling only the inflow $u_1$, i.e., for what values of $u_3^* \in (0, 0.5)$ do we have $R = \{1\}$ ?".

We checked numerically inequalities (3.11) by computing the values of $\mu_i$ from (3.6) for given values of $u_3^*$. It was found that inequalities (3.11) hold for $R = \{1\}$, provided that $u_3^* < 0.1$. Therefore, we conclude that the uncongested equilibrium point can be globally exponentially stabilized by controlling only the inflow $u_1$ for $u_3^* < 0.1$. The answer may be conservative, since the estimation of the constant $C > 0$ that satisfies (3.2) is conservative. ◁

**Example 4.2:** Consider a freeway model of the form (2.11), (2.12), (2.13), (2.8) with $n = 5$ cells. Each cell has the same critical density $\delta_i = 55$ ($i = 1,...,5$) and the same jam density $a_i = 170$ ($i = 1,...,5$). The considered freeway stretch has no intermediate on/off-ramps (i.e., $u_i(t) \equiv u_i^* = 0$ for $i = 2, 3, 4, 5$ and $p_i = 0$ for $i = 1,...,4$). Thus, the only control possibility is the inflow $u_1$ of the first cell (see Figure 3). We also suppose that the cell flow capacities are $q_i = 25$ for $i = 1, 2, 3, 4$ and $q_5 = 20$, i.e. the last cell has 20% lower flow capacity (e.g. due to grade or curvature or tunnel or bridge etc.) than the first four cells and is therefore a potential bottleneck for the freeway.

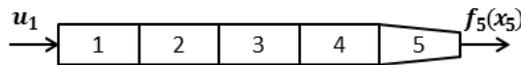

**Figure 3:** Freeway stretch.

Figure 4, depicts the triangular fundamental diagrams for the above model. The blue line corresponds to the demand function, while the red line corresponds to the supply function. More precisely, the demand part for every cell is given by the following functions:

$$f_i(z) = \begin{cases} (5/11)z & z \in [0,55] \\ (25/115)(170-z) & z \in (55, 87.2] \\ 18 & z \in (87.2, 170] \end{cases}$$

$$(i=1,...,4), f_5(z) = \begin{cases} (4/11)z & z \in [0,55] \\ (20/115)(170-z) & z \in (55, 72.25] \\ 17 & z \in (72.25, 170] \end{cases}$$

Notice that, the capacity drop phenomenon has been taken into account by considering a partly decreasing demand function for over-critical densities $x_i \in (55, 170]$.

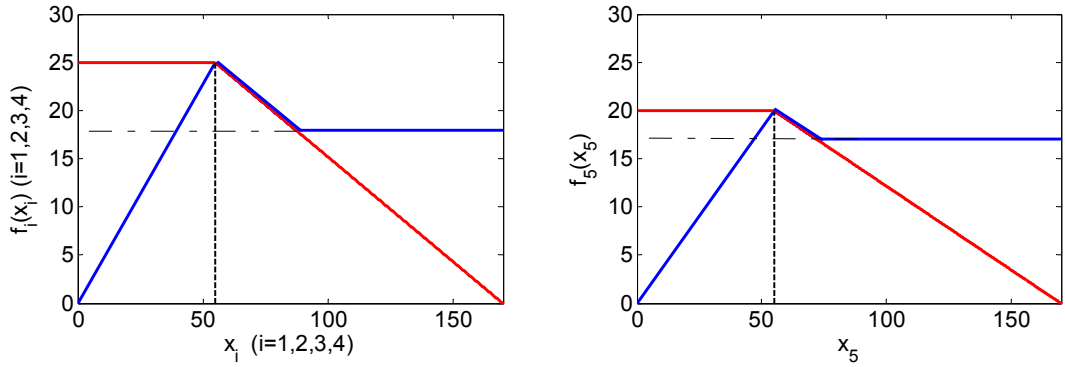

**Figure 4:** Fundamental diagram of every cell.

For readers who are accustomed to the traditional units of veh/h for flows and veh/km for densities, the example model may be viewed to reflect a freeway stretch with 3-lane cells with equal lengths of 500 m; with a time step of 15 s. With these settings, the critical density of $\delta_i = 55$ corresponds to 36.7 veh/km/lane; while the jam density of $a_i = 170$ corresponds to 113.3 veh/km/lane; and the cell flow capacities of $q_{1,2,3,4} = 25$ and $q_5 = 20$ correspond to 2000 veh/h/lane and 1600 veh/h/lane, respectively.

Assumption (H) holds with $\delta_i = \tilde{\delta}_i = 55$ ($i=1,...,5$), $L_i = 6/11$ ($i=1,...,4$), $L_5 = 7/11$. The uncongested equilibrium point $x_i^* = 11u_1^*/5$ ($i=1,...,4$), $x_5^* = 11u_1^*/4$ exists for $u_1^* < 20$. Simulations showed that the open-loop system converges to an uncongested equilibrium point for main inflow $u_1^*$ less than 17. For higher values of the main inflow, the uncongested equilibrium point is not globally exponentially stable due to the existence of additional (congested) equilibria. This is shown in Figure 5, where the solution of the open-loop system, with constant inflow $u_1^* = 19.99$ is attracted by the congested equilibrium $[91.8, 91.8, 91.8, 91.8, 72.25]$. The components of the uncongested equilibrium for $u_1^* = 19.99$ are $x_i^* = 43.978$ for $i=1,...,4$ and $x_5^* = 54.9725$. Therefore, if the objective is the operation of the freeway with large flows, then a control strategy will be needed.

We next notice that property (C2) holds with $\lambda_i = 6/11$, $G_i = 5/11$ ($i = 1,...,4$), $\lambda_5 = 7/11$ and $G_5 = 4/11$. Therefore, we are in a position to achieve global exponential stabilization of the uncongested equilibrium point for model by using Theorem 2.1. Indeed, Theorem 2.1 guarantees that for every $\sigma \in (0,1]$ there exists a constant $b_1 \in (0, u_1^*)$ and a constant $\gamma > 0$ such that, the feedback law $k : (0,10]^5 \to \Re_+$ defined by:

$$u_1 = \max\left(u_1^* - \gamma \sum_{i=1}^{5} \sigma^i \max\left(0, x_i - x_i^*\right), b_1\right) \qquad (4.1)$$

achieves robust global exponential stabilization of the uncongested equilibrium point $x^* = (x_1^*,..., x_5^*) \in (0,55) \times ... (0,55)$ for the closed-loop system.

We selected $u_1^* = 19.99$, which is very close to $20$, the capacity flow of cell 5. The value of the constant $b_1 \in (0, u_1^*)$ was chosen to be 0.2; this is a rather low minimum flow value in practice, but allows us here to study the dynamic properties of the regulators in a broader feasible control area. Various values of the constants $\sigma \in (0,1]$ and $\gamma > 0$ were tested by performing a simulation study with respect to various initial conditions. Low values for $\sigma \in (0,1]$ require large values for $\gamma > 0$ in order to have global exponential stability for the closed-loop system. Moreover, in order to evaluate the performance of the controller, we used as a performance criterion the total number of Vehicles Exiting the Freeway (*VEF*) on the interval $[0,T]$, i.e.,

$$VEF_T = \sum_{t=0}^{T} f_5(x_5(t)). \qquad (4.2)$$

Notice that the freeway performs best (and total delays are minimised) if *VEF* is maximized; the maximum theoretical value for *VEF* is $20(T+1)$, which is achieved if cell 5 is operating at capacity flow ($q_5 = 20$) at all times. For $T = 200$, the maximum theoretical value of *VEF* is 4020.

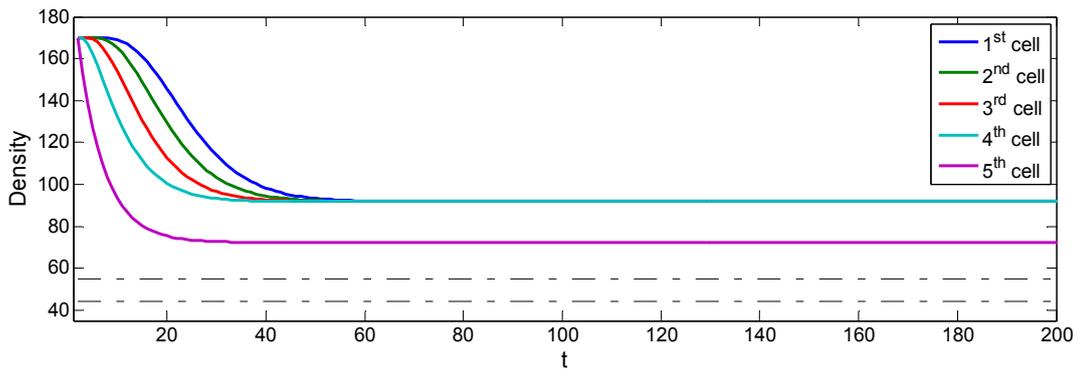

**Figure 5:** Open-loop system convergence (dashed lines correspond to the uncongested equilibrium point for the inflow $u_1^* = 19.99$).

All following tests of the proposed regulator (4.1) were conducted with the same values $\sigma = 0.7$ and $\gamma = 0.6$. The responses of the densities of every cell for the closed-loop system with the proposed feedback regulator (4.1) and initial condition $x_0 = [60, 57, 58, 60, 62]$ are shown in Figure 6(a). Notice that all initial cell density values are slightly overcritical (slightly congested).

For this case, we had $VEF_{200} = 3979.8$. The feedback regulator is seen to respond very satisfactorily in this test and achieves an accordingly high performance.

A detailed comparison of the proposed feedback regulator (4.1) is one of the very few was made with the Random Located Bottleneck (RLB) PI regulator, which was proposed in [32] and comparable feedback regulators that has been employed in field operations [27]. The RLB PI regulator for the present system is implemented as follows:

$$v_i(t) = \min\left(u_{\max}, \min\left(q_1, c_1(a_1 - x_1(t-1)), u_1(t-1)\right) + c, \max\left(u_{\min}, v_i(t-1) - K_p(x_i(t) - x_i(t-1)) + K_I(\delta_i - x_i(t))\right)\right) \quad (4.3)$$

for $i = 1, \ldots, 5$,

$$v_i^{sm}(t) = e v_i(t) + (1-e) v_i^{sm}(t-1), \text{ for } i = 1, \ldots, 5 \quad (4.4)$$

$$j(t) = \min\left\{l \in \{1, 2, 3, 4, 5\} : v_l^{sm}(t) = \min_{i=1,\ldots,5}\left(v_i^{sm}(t)\right)\right\} \quad (4.5)$$

$$u_1(t) = v_{j(t)}(t). \quad (4.6)$$

where $c, K_p, K_I > 0$ and $e \in (0,1)$ are constant parameters. Essentially, (4.3) reflects the parallel (independent) operation of five bounded PI-type regulators, one for each cell; while (4.4) performs an exponential smoothing of the respective obtained inflows. Eventually, the smoothed inflow values are compared in (4.5) in order to pick the currently most conservative regulator; whose (unsmoothed) inflow is finally actually activated as a control input in (4.8), see [32] for the background and detailed reasoning for this approach. The parameters for the RLB PI regulators are set (as proposed in [32] - with the suitable transformation in the current units) to be $K_p = 5/18$, $K_I = 1/90$, while $c = 4$, $e = 0.5$, $u_{\min} = 0.2$ and $u_{\max} = 25$. Notice that all PI regulators were given the same gain values for simplicity and convenience, as suggested in [32]. In all reported tests, the initial condition for the RLB PI regulator was $v_i(-1) = v_i^{sm}(-1) = u_1(-1) = 20$ for $i = 1, \ldots, 5$, and $x(-1) = x(0) = x_0$, where $x_0$ is the vector of the initial values for the densities of every cell.

When applied to the same initial condition $x_0 = [60, 57, 58, 60, 62]$, the RLB PI regulator (Figure 6(b)), led to slower convergence compared with the proposed regulator (4.1). This is also reflected in the computed value of $VEF_{200} = 3785.9$ for RLB PI regulator. In general, conducting a simulation study with various levels of initial conditions, the proposed regulator (4.1) exhibited faster performance than the RLB PI regulator. For example, Figure 7 shows the evolution of the Euclidean norm $|x(t) - x^*|$ for the closed-loop system with the proposed feedback regulator (4.1) (blue curve) and for the closed-loop system with the RLB PI regulator (4.3), (4.4), (4.5), (4.6) (red curve), when starting from the initial condition $[170, 170, 170, 170, 170]$ reflecting a fully congested original state (as in Figure 5 as well). It is again clear that the proposed feedback regulator (4.1) achieves faster convergence and higher performance of $VEF_{200} = 3845.2$, while $VEF_{200} = 3007.8$ resulted for the RLB PI regulator.

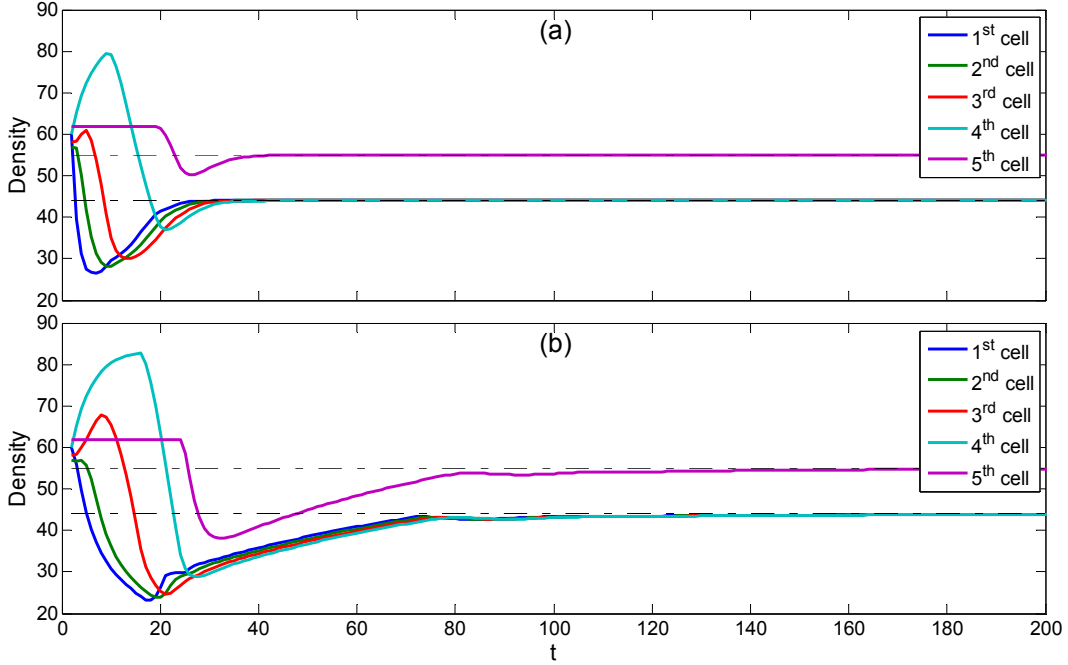

**Figure 6:** The responses of the densities of every cell for the closed-loop system (4.1) with initial condition $x_0 = [60, 57, 58, 60, 62]$ using: (a) the proposed feedback regulator (4.1); and (b) the RLB PI regulator.

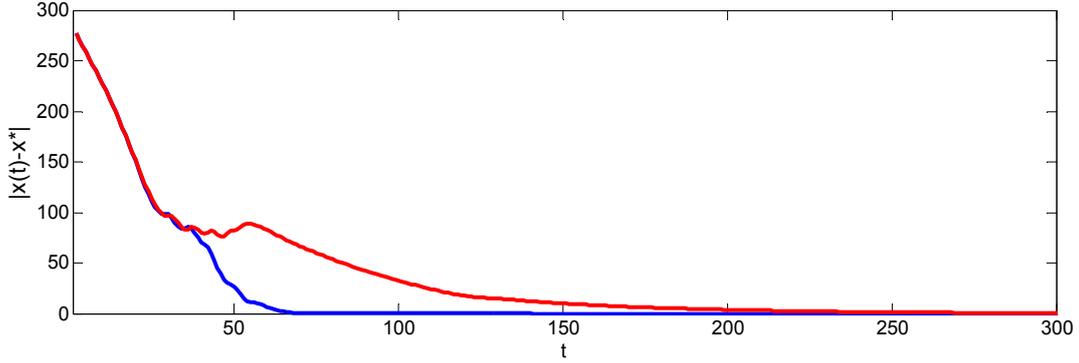

**Figure 7:** The evolution of the Euclidean norm for the closed-loop system (4.1) and initial condition $x_0 = [170, 170, 170, 170, 170]$ for two cases, for the proposed feedback regulator (4.1) (blue curve) and for the RLB PI regulator (red curve).

We next investigated the robustness of the proposed feedback regulator with respect to measurement errors. The applied formula for the measurements was:

$$\tilde{x}(t) = P(x(t) + Ae(t)) \tag{4.7}$$

where $P$ is the projection operator on the closure of $S$, $e(t)$ is a normalized vector, and $A \geq 0$ is the magnitude of the measurement error. In this case, the feedback law (4.1) was implemented based on the state measurement $\tilde{x}(t)$ given by (4.7), i.e.,

$$u_1(t) = \max\left( u_1^* - \gamma \sum_{i=1}^{5} \sigma^i \max\left(0, \tilde{x}_i(t) - x_i^*\right), b_1 \right). \tag{4.8}$$

For comparison purposes, we also present the performance of the RLB PI regulator for the same system, under the same measurement errors. In this case, equation (4.3) is replaced by the equation

$$v_i(t) = \min\left(u_{\max}, \min(q_1, c_1(a_1 - x_1(t-1)), u_1(t-1)) + c, \max\left(u_{\min}, v_i(t-1) - K_p(\tilde{x}_i(t) - \tilde{x}_i(t-1)) + K_I(\delta_i - \tilde{x}_i(t))\right)\right)$$
(4.9)

for $i = 1,...,5$, where the state measurement $\tilde{x}(t)$ is given by (4.7).

Figure 8 shows the responses of the densities of every cell for two cases: (a) for the closed-loop system with the proposed feedback regulator (4.8); and (b) for the closed-loop system with the RLB PI regulator (4.9), (4.4), (4.5), (4.6); where the state measurement in both cases is given by (4.7) with $A = 10$, $e(t) = \dfrac{\cos(\omega t)}{\sqrt{5}}(1,1,...,1)$, $\omega = \pi$. The initial condition was the uncongested equilibrium point.

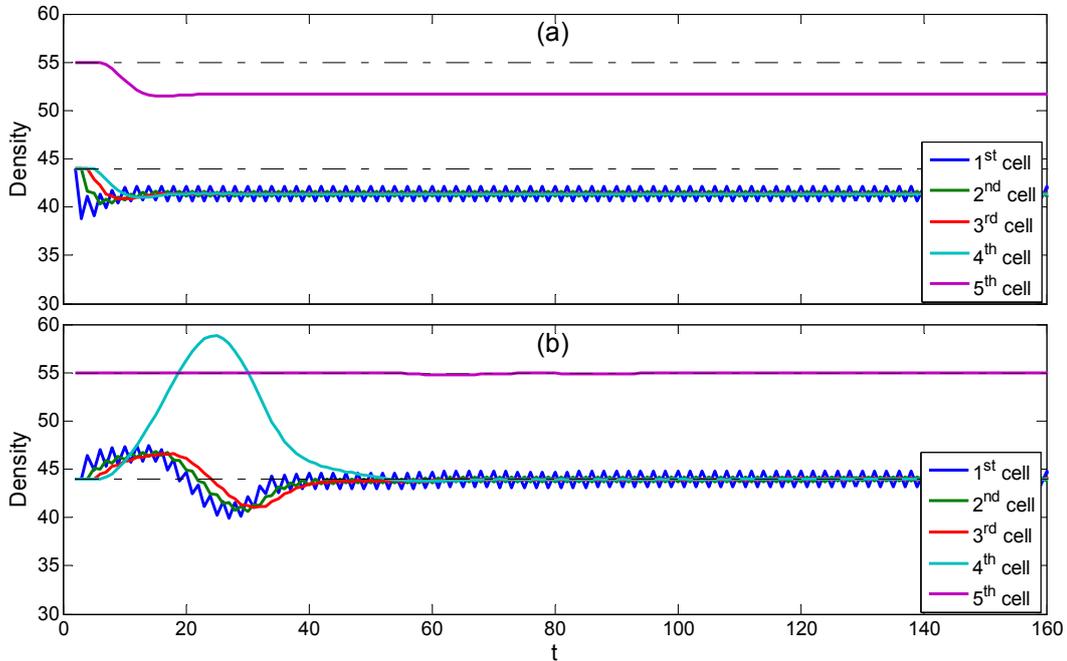

**Figure 8:** The responses of the densities of every cell for the closed-loop system (4.1): (a) with the proposed feedback regulator (4.8); and (b) with the RLB PI regulator (4.9), (4.4), (4.5), (4.6). In both cases the state measurement is given by (4.7) with $A = 10$, $e(t) = \left(\cos(\omega t)/\sqrt{5}\right)(1,1,...,1)$, $\omega = \pi$; initial condition is the uncongested equilibrium point.

In this test, the RLB PI regulator is less sensitive to measurement errors than the proposed feedback regulator (4.8), the latter producing a visible offset (Figure 8), in particular also for cell 5, which reduces accordingly the stationary outflow. This is also reflected in the computed values of $VEF_{200} = 3789$ for the proposed feedback regulator (4.9) (which is 6% less than the maximum value of $VEF_{200}$) and $VEF_{200} = 4016.8$ for the RLB PI regulator (which is 0.8% less than the maximum value of $VEF_{200}$) due to the measurement error. The ultimate mean values of the states are much closer to the equilibrium values for the RLB PI regulator than for the proposed feedback regulator (4.8), indicating that the RLB PI regulator achieves a much smaller mean offset in this case. It should be noted at this point that various frequencies $\omega$ were tested for measurement errors. While Figure 8 is typical for medium and high frequencies (the RLB PI regulator achieves a smaller mean offset than the proposed feedback regulator (4.8)), the results indicate higher sensitivity of the RLB PI regulator with respect to measurement errors at low frequencies (Figure 9). For low frequency measurement errors, the proposed feedback regulator (4.8) achieves a smaller mean offset than the RLB PI regulator, as shown in Figure 9. Figure 9 shows the responses of the densities of every cell for two cases: (a) for the closed-loop system with the proposed feedback regulator (4.8), and (b) for the closed-loop system with the RLB PI regulator (4.9), (4.4), (4.5), (4.6), where the state measurement in both cases is given by (4.7) with $A=10$, $e(t)=\left(\cos(\omega t)/\sqrt{5}\right)(1,1,...,1)$, $\omega=0.1$. The initial condition is the uncongested equilibrium point.

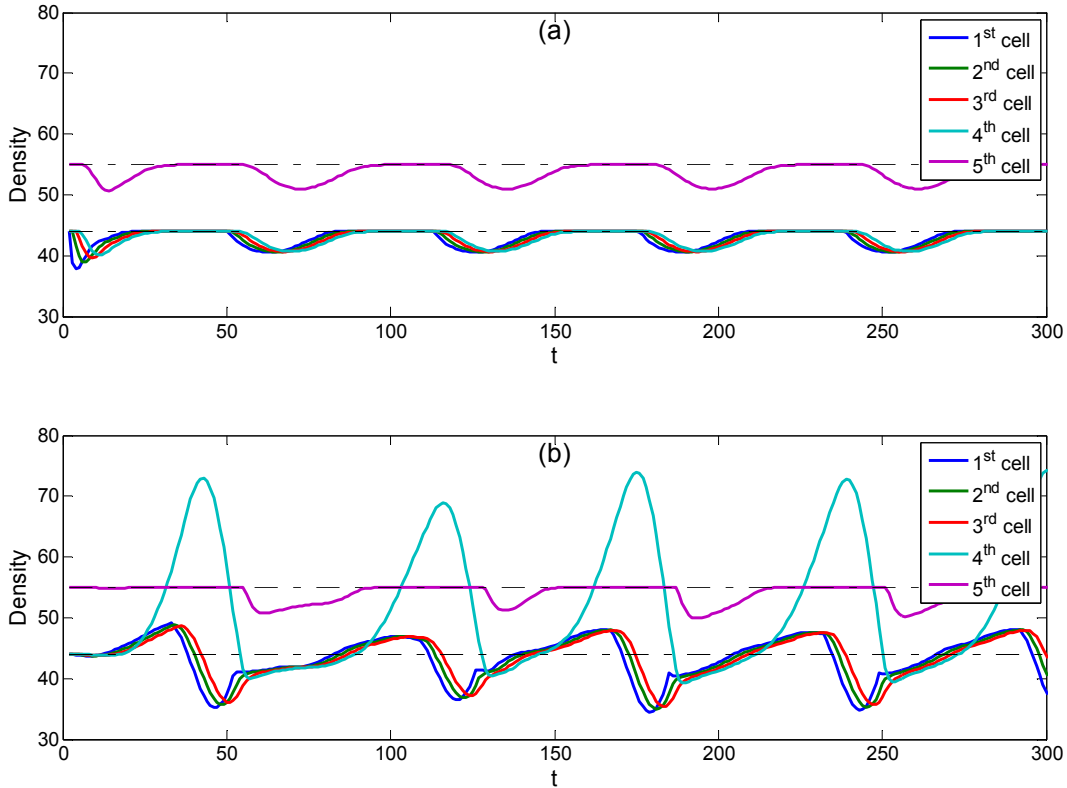

**Figure 9:** The responses of the densities of every cell for the closed-loop system (4.1): (a) with the proposed feedback regulator (4.8); and (b) with the RLB PI regulator (4.9), (4.4), (4.5), (4.6). In both cases the state measurement is given by (4.7) with $A=10$, $e(t)=\left(\cos(\omega t)/\sqrt{5}\right)(1,1,...,1)$, $\omega=0.1$; initial condition is the uncongested equilibrium point.

The conclusions of this simulation study are:

- The proposed feedback regulator (4.1) can achieve a faster convergence of the state to the equilibrium compared to the RLB PI regulator in the absence of measurement errors.

- The proposed feedback regulator (4.1) is quite robust to measurement errors. However, it is more sensitive to measurement errors with high frequency than the RLB PI regulator; but it is less sensitive to low-frequency measurement errors than the RLB PI regulator. Intended future extensions are expected to improve the properties of the proposed feedback regulator in this respect, as well as in cases of modelling errors or persisting disturbances.

# 5. Conclusions

This work provided a rigorous methodology for the construction of a parameterized family of explicit feedback laws that guarantee the robust global exponential stability of the uncongested equilibrium point for general nonlinear and uncertain discrete-time freeway models. The construction of the global exponential feedback stabilizer was based on the CLF (control Lyapunov function) approach as well as on certain important properties of freeway models.

Simulation-based comparisons were made with existing feedback laws, which were proposed in the literature and have been in practical use. More specifically, we compared the performance and some robustness properties of the closed-loop system under the effect of the proposed feedback law and under the effect of the Random Located Bottleneck (RLB) PI regulator [32]. In most cases, it was found that the performance and the robustness properties guaranteed by the implementation of the proposed feedback law were good and comparable to or better than the performance and the robustness properties induced by the RLB PI regulator.

Ongoing and future research addresses robustness issues in a rigorous way: the knowledge of a Lyapunov function for the closed-loop system can be exploited to this purpose, and explicit formulas for the gains of various inputs (measurement or modelling errors) can be derived. Also, the estimation of the gains of various inputs can allow the study and control of interconnected freeways (traffic networks). Finally, the present approach does not consider the impact of inflow control on upstream traffic flow conditions (e.g. queue forming at on-ramps); future extensions will address these issues appropriately.

## Acknowledgments

The research leading to these results has received funding from the European Research Council under the European Union's Seventh Framework Programme (FP/2007-2013) / ERC Grant Agreement n. [321132], project TRAMAN21.

## Appendix

**Proof of (C4):** We prove the following claim:

**(Claim):** For all $m = 1,...,n-1$ there exists a constant $C_m > 0$ such that the following inequality holds for all $x \in S := (0, a_1] \times ... \times (0, a_n]$, $u \in U = (0, +\infty) \times [0, r_2] \times ... \times [0, r_n]$, $d = (d_2, ..., d_n) \in [0,1]^{n-1}$:

$$\sum_{i=m}^{n} (1 + p_i(n-i)) s_{i+1} f_i(x_i) \geq C_m \sum_{i=m}^{n} (n+1-i) x_i. \quad (A.1)$$

Property (C4) is a direct consequence of the above claim.

First we prove the claim for $m = n-1$. Define $l_n = \min\left(1, \frac{c_n a_n - r_n}{2(1 - p_{n-1}) f_{n-1}(\delta_{n-1})}\right)$. Indeed, using Property (C3), we get for all $(x_{n-1}, x_n) \in (0, a_{n-1}] \times (0, a_n]$, $u \in U = (0, +\infty) \times [0, r_2] \times ... \times [0, r_n]$, $d = (d_2, ..., d_n) \in [0,1]^{n-1}$ with $s_n \geq l_n$:

$$f_n(x_n) + (1 + p_{n-1}) s_n f_{n-1}(x_{n-1}) \geq$$
$$l_n (1 + p_{n-1}) \theta_{n-1} x_{n-1} + \theta_n x_n \geq (x_n + 2 x_{n-1}) \min\left(\theta_n, \frac{1 + p_{n-1}}{2} l_n \theta_{n-1}\right) \quad (A.2)$$

On the other hand, for all $(x_{n-1}, x_n) \in (0, a_{n-1}] \times (0, a_n]$, $u \in U = (0, +\infty) \times [0, r_2] \times ... \times [0, r_n]$, $d = (d_2, ..., d_n) \in [0,1]^{n-1}$ with $s_n < l_n$, it follows from (2.8) that $(1 - p_{n-1}) f_{n-1}(x_{n-1}) + u_n > \min(q_n, c_n(a_n - x_n))$. We distinguish the cases:

- $\min(q_n, c_n(a_n - x_n)) = c_n(a_n - x_n)$. In this case, we have $a_n - c_n^{-1} q_n \leq x_n$. Assumption (H) in conjunction with the fact that $u_n \leq r_n < c_n a_n$ implies that $l_n (1 - p_{n-1}) f_{n-1}(\delta_{n-1}) + r_n > c_n(a_n - x_n)$, which by virtue of definition $l_n = \min\left(1, \frac{c_n a_n - r_n}{2(1 - p_{n-1}) f_{n-1}(\delta_{n-1})}\right)$, gives $x_n > \frac{1}{2}\left(a_n - c_n^{-1} r_n\right)$. Using Property (C3), we get:

$$f_n(x_n) + (1 + p_{n-1}) s_n f_{n-1}(x_{n-1}) \geq \theta_n x_n \geq \theta_n \frac{1}{2}\left(a_n - c_n^{-1} r_n\right) \geq \theta_n \frac{1}{2}\left(a_n - c_n^{-1} r_n\right) \frac{x_n + 2 x_{n-1}}{a_n + 2 a_{n-1}}$$

- $\min(q_n, c_n(a_n - x_n)) = q_n$. In this case, we obtain from (2.8) and the fact that $u_n \leq r_n < q_n$ the inequality $s_n \geq \dfrac{q_n - r_n}{(1 - p_{n-1}) f_{n-1}(x_{n-1})}$. Using Property (C3), we get:

$$f_n(x_n) + (1 + p_{n-1}) s_n f_{n-1}(x_{n-1}) \geq \theta_n x_n + (1 + p_{n-1}) \frac{q_n - r_n}{1 - p_n}$$

$$\geq \theta_n x_n + (1 + p_{n-1}) \frac{q_n - r_n}{1 - p_n} \frac{2 x_{n-1}}{2 a_{n-1}} \geq (x_n + 2 x_{n-1}) \min\left(\theta_n, \frac{q_n - r_n}{1 - p_n} \frac{1 + p_{n-1}}{2 a_{n-1}}\right)$$

It follows from (A.2) and the above inequalities that (A.1) holds with

$$C_{n-1} := \min\left(\theta_n, \frac{1 + p_{n-1}}{2} l_n \theta_{n-1}, \frac{\theta_n (a_n - c_n^{-1} r_n)}{2(a_n + 2 a_{n-1})}, \frac{1 + p_{n-1}}{1 - p_{n-1}} \frac{q_n - r_n}{2 a_{n-1}}\right).$$

Next, we suppose that the claim holds for $m = k \in \{2, \ldots, n-1\}$ and we show that the claim holds for $m = k - 1$. Define $l_k = \min\left(1, \dfrac{c_k a_k - r_k}{2(1 - p_{k-1}) f_{k-1}(\delta_{k-1})}\right)$. Using (A.1) for $m = k$ and Property (C3), we obtain for all $x \in S := (0, a_1] \times \ldots \times (0, a_n]$, $u \in U = (0, +\infty) \times [0, r_2] \times \ldots \times [0, r_n]$, $d = (d_2, \ldots, d_n) \in [0, 1]^{n-1}$ with $s_k \geq l_k$:

$$\sum_{i=k-1}^{n} (1 + p_i(n - i)) s_{i+1} f_i(x_i)$$

$$= \sum_{i=k}^{n} (1 + p_i(n - i)) s_{i+1} f_i(x_i) + (1 + p_{k-1}(n + 1 - k)) s_k f_{k-1}(x_{k-1})$$

$$\geq C_k \sum_{i=k}^{n} (n + 1 - i) x_i + (1 + p_{k-1}(n + 1 - k)) l_k \theta_{k-1} x_{k-1} \geq$$

$$\geq \min\left(C_k, \frac{1 + (n + 1 - k) p_{k-1}}{n + 2 - k} l_k \theta_{k-1}\right) \sum_{i=k-1}^{n} (n + 1 - i) x_i$$

(A.3)

On the other hand, for all $x \in S := (0, a_1] \times \ldots \times (0, a_n]$, $u \in U = (0, +\infty) \times [0, r_2] \times \ldots \times [0, r_n]$, $d = (d_2, \ldots, d_n) \in [0, 1]^{n-1}$ with $s_k < l_k$, it follows from (2.8) that $l_k (1 - p_{k-1}) f_{k-1}(x_{k-1}) + u_k > \min(q_k, c_k(a_k - x_k))$. We distinguish the cases:

- $\min(q_k, c_k(a_k - x_k)) = c_k(a_k - x_k)$. Assumption (H) in conjunction with the fact that $u_k \leq r_k < c_k a_k$ implies that $l_k(1 - p_{k-1}) f_{k-1}(\delta_{k-1}) + r_k > c_k(a_k - x_k)$, which by virtue of definition $l_k = \min\left(1, \dfrac{c_k a_k - r_k}{2(1 - p_{k-1}) f_{k-1}(\delta_{k-1})}\right)$, gives $x_k > \dfrac{1}{2}(a_k - c_k^{-1} r_k)$. Using (A.1) for $m = k$ and Property (C3), we get:

$$\sum_{i=k-1}^{n} \left(1 + p_i(n-i)\right) s_{i+1} f_i(x_i) \geq \sum_{i=k}^{n} \left(1 + p_i(n-i)\right) s_{i+1} f_i(x_i)$$

$$\geq C_k \sum_{i=k}^{n} (n+1-i) x_i = C_k \sum_{i=k+1}^{n} (n+1-i) x_i + C_k (n+1-k) x_k$$

$$\geq C_k \sum_{i=k+1}^{n} (n+1-i) x_i + C_k (n+1-k) \frac{1}{2} \left(a_k - c_k^{-1} r_k\right)$$

$$\geq C_k \sum_{i=k+1}^{n} (n+1-i) x_i + C_k (n+1-k) \frac{1}{2} \left(a_k - c_k^{-1} r_k\right) \frac{(n+1-k) x_k + (n+2-k) x_{k-1}}{(n+1-k) a_k + (n+2-k) a_{k-1}}$$

$$\geq C_k \min\left(1, \frac{(n+1-k)\left(a_k - c_k^{-1} r_k\right)}{2(n+1-k) a_k + 2(n+2-k) a_{k-1}}\right) \sum_{i=k-1}^{n} (n+1-i) x_i$$

- $\min(q_k, c_k(a_k - x_k)) = q_k$. In this case, we obtain from (2.8) and the fact that $u_k \leq r_k < q_k$ the inequality $s_k \geq \dfrac{q_k - r_k}{(1 - p_{k-1}) f_{k-1}(x_{n-1})}$. Consequently, we get:

$$\sum_{i=k-1}^{n} \left(1 + p_i(n-i)\right) s_{i+1} f_i(x_i)$$

$$= \sum_{i=k}^{n} \left(1 + p_i(n-i)\right) s_{i+1} f_i(x_i) + \left(1 + p_{k-1}(n+1-k)\right) s_k f_{k-1}(x_{k-1})$$

$$\geq C_k \sum_{i=k}^{n} (n+1-i) x_i + \left(1 + p_{k-1}(n+1-k)\right) \frac{q_k - r_k}{1 - p_{k-1}}$$

$$\geq C_k \sum_{i=k}^{n} (n+1-i) x_i + \left(1 + p_{k-1}(n+1-k)\right) \frac{q_k - r_k}{1 - p_{k-1}} \frac{(n+2-k) x_{k-1}}{(n+2-k) a_{k-1}}$$

$$\geq \min\left(C_k, \frac{q_k - r_k}{1 - p_{k-1}} \frac{1 + p_{k-1}(n+1-k)}{(n+2-k) a_{k-1}}\right) \sum_{i=k-1}^{n} (n+1-i) x_i$$

It follows from (A.3) and the above inequalities that the claim holds for $m = k-1$ with:

$$C_{k-1} = \min\left(C_k, \frac{1 + p_{k-1}(n+1-k)}{n+2-k} l_k \theta_{k-1}, \frac{(n+1-k) C_k \left(a_k - c_k^{-1} r_k\right)}{2(n+1-k) a_k + 2(n+2-k) a_{k-1}}, \frac{q_k - r_k}{1 - p_{k-1}} \frac{1 + p_{k-1}(n+1-k)}{(n+2-k) a_{k-1}}\right). \triangleleft$$

**Proof of (C5):** The following equations hold for all $x \in S := (0, a_1] \times \ldots (0, a_n]$, $u = (u_1, \ldots, u_n)' \in (0, +\infty) \times \Re_+^{n-1}$, $d = (d_2, \ldots, d_n) \in [0,1]^{n-1}$ with $p_n = 1 = s_{n+1}$ and are direct consequences of (2.11), (2.12), (2.13) and definitions $I_j(x) := \sum_{i=1}^{j} x_i$ for $j = 1, \ldots, n$:

$$I_j(x^+) = I_j(x) + \sum_{i=1}^{j} w_i u_i - \sum_{i=1}^{j-1} s_{i+1} p_i f_i(x_i) - s_{j+1} f_j(x_j), \text{ for } j = 2, \ldots, n-1 \quad \text{(A.4)}$$

$$I_n(x^+) = I_n(x) + \sum_{i=1}^{n} w_i u_i - \sum_{i=1}^{n-1} s_{i+1} p_i f_i(x_i) - f_n(x_n). \tag{A.5}$$

Equality (3.4) is a consequence of (2.11), (A.4), (A.5) and definitions $p_n = 1 = s_{n+1}$.

Combining (3.4) and (3.2), we get:

$$\sum_{i=1}^{n} I_i(x^+) \leq \sum_{i=1}^{n} I_i(x) - C \sum_{i=1}^{n} (n+1-i) x_i + \sum_{i=1}^{n} (n+1-i) w_i u_i, \text{ for all } (x,u,d) \in S \times U \times [0,1]^{n-1}. \tag{A.6}$$

Since $w_i \in [0,1]$ ($i = 1,...,n$) and $\sum_{i=1}^{n} I_i(x) = \sum_{i=1}^{n} (n+1-i) x_i$, it follows from (A.6) that (3.5) holds. ◁

**Proof of the Claim 1 made in the proof of Theorem 2.1:** We distinguish two cases.

Case 1: $x \in \Omega = (0, \mu_1] \times ... \times (0, \mu_n]$, $d \in [0,1]^{n-1}$.

Definition (2.20) and equations (3.4), (3.8), (3.9) with $u_i = k_i(x) \leq u_i^*$ give:

$$V(x^+) = \sigma \left| x_1 - f_1(x_1) + u_1 - x_1^* \right| + \sum_{i=2}^{n} \sigma^i \left| x_i - f_i(x_i) + (1 - p_{i-1}) f_{i-1}(x_{i-1}) + u_i - x_i^* \right|$$
$$+ K \max\left( 0, \sum_{i=1}^{n} I_i(x) - \sum_{i=1}^{n} (1 + (n-i) p_i) f_i(x_i) + \sum_{i=1}^{n} (n+1-i) u_i - P(x^+) \right) + A\Xi(x^+) \tag{A.7}$$

with $p_n = 1$. Using (3.13), property (C2), the fact that $\mu_i \leq \tilde{\delta}_i$ for $i = 1,...,n$ (a consequence of (3.6) and (*)) and definition (3.10), we get from (A.7):

$$V(x^+) \leq L \sum_{i=1}^{n} \sigma^i \left| x_i - x_i^* \right| + \sum_{i=1}^{n} \sigma^i \left| u_i - u_i^* \right| + LA\Xi(x)$$
$$+ K \max\left( 0, \sum_{i=1}^{n} I_i(x) - \sum_{i=1}^{n} (1 + (n-i) p_i) f_i(x_i) + \sum_{i=1}^{n} (n+1-i) u_i - P(x^+) \right) \tag{A.8}$$

It follows from the combination of (2.18) and inequality (A.8) that the following inequality holds for all $x \in \Omega$:

$$V(x^+) \leq L \sum_{i=1}^{n} \sigma^i \left| x_i - x_i^* \right| + \sum_{i \in R} \sigma^i \min\left( \gamma_i \Xi(x), u_i^* - b_i \right) + LA\Xi(x)$$
$$+ K \max\left( 0, \sum_{i=1}^{n} I_i(x) - \sum_{i=1}^{n} (1 + (n-i) p_i) f_i(x_i) + \sum_{i=1}^{n} (n+1-i) u_i - P(x^+) \right) \tag{A.9}$$

Inequality (3.2) and equations (3.8) imply that:

$$\sum_{i=1}^{n}(1+(n-i)p_i)f_i(x_i) \geq C\sum_{i=1}^{n}I_i(x). \tag{A.10}$$

Using (A.8) and (A.10), we get:

$$\begin{aligned}V(x^+) &\leq L\sum_{i=1}^{n}\sigma^i|x_i-x_i^*| + \sum_{i\in R}\sigma^i \min\left(\gamma_i\Xi(x), u_i^*-b_i\right) + LA\Xi(x) \\ &+ K\max\left(0, (1-C)\sum_{i=1}^{n}I_i(x) + \sum_{i=1}^{n}(n+1-i)u_i - P(x^+)\right)\end{aligned} \tag{A.11}$$

We next distinguish two cases:

Case 1(i): $\Xi(x) \leq \tau$. In this case we have $\gamma_i\Xi(x) \leq u_i^* - b_i$ for all $i \in R$. Since $\Xi(x) \leq h$ (a consequence of $\tau < \tau^* \leq h$), we get from (3.13) and definition (2.18) that $u_i = k_i(x) = u_i^* - \gamma_i\Xi(x) \geq b_i$ for all $i \in R$ and $\min(h, \Xi(x^+)) \leq L\min(h, \Xi(x))$. Using the definitions $I_j(x) := \sum_{i=1}^{j}x_i$ for $j = 1,\ldots,n$, $P(x) := Q - \theta\min(h, \Xi(x))$ and the facts

- $Q \geq (1-C)\sum_{i=1}^{n}I_i(x^*) + (1-C)h\max_{i=1,\ldots,n}\left((n+1-i)\sigma^{-i}\right) + \sum_{i=1}^{n}(n+1-i)u_i^*$

- $\sum_{i\in R}(n+1-i)\gamma_i = \tau^{-1}\sum_{i\in R}(n+1-i)(u_i^*-b_i) \geq (\tau^*)^{-1}\sum_{i\in R}(n+1-i)(u_i^*-b_i) \geq \theta L$ (a consequence of $\tau \leq \tau^* \leq (\theta L)^{-1}\sum_{i\in R}(n+1-i)(u_i^*-b_i)$),

- $\sum_{i=1}^{n}(n+1-i)(x_i-x_i^*) \leq \sum_{i=1}^{n}\left((n+1-i)\sigma^{-i}\right)\sigma^i\max(0, x_i-x_i^*) \leq \max_{i=1,\ldots,n}\left((n+1-i)\sigma^{-i}\right)\Xi(x)$ for $i = 1,\ldots,n$ (a consequence of definition (2.19)),

we get:

$$Q \geq (1-C)\sum_{i=1}^{n}I_i(x^*) + (1-C)h\max_{i=1,\ldots,n}\left((n+1-i)\sigma^{-i}\right) + \sum_{i=1}^{n}(n+1-i)u_i^* \Rightarrow$$

$$Q \geq (1-C)\sum_{i=1}^{n}I_i(x^*) + (1-C)\Xi(x)\max_{i=1,\ldots,n}\left((n+1-i)\sigma^{-i}\right) + \sum_{i=1}^{n}(n+1-i)u_i^* \Rightarrow$$

$$Q \geq \theta L\Xi(x) - \sum_{i\in R}(n+1-i)\gamma_i\Xi(x) + (1-C)\sum_{i=1}^{n}I_i(x^*) + (1-C)\sum_{i=1}^{n}(n+1-i)(x_i-x_i^*) + \sum_{i=1}^{n}(n+1-i)u_i^* \Rightarrow$$

$$Q \geq \theta\min(h, \Xi(x^+)) + (1-C)\sum_{i=1}^{n}I_i(x) + \sum_{i\in R}(n+1-i)(u_i^*-\gamma_i\Xi(x)) + \sum_{i\notin R}(n+1-i)u_i^* \Rightarrow$$

$$0 \geq \sum_{i=1}^{n}(n+1-i)u_i + (1-C)\sum_{i=1}^{n}I_i(x) - P(x^+)$$

Combining (A.11) with the above inequality, we obtain:

$$V(x^+) \leq L\sum_{i=1}^{n} \sigma^i |x_i - x_i^*| + \sum_{i \in R} \sigma^i \gamma_i \Xi(x) + LA\Xi(x). \quad (A.12)$$

It follows from (A.12) and the fact that $A \geq (1-L)^{-1} \sum_{i \in R} \sigma^i \gamma_i$ that (3.18) holds when $\Xi(x) \leq \tau$.

Case 1(ii): $\Xi(x) > \tau$. In this case $\gamma_i \Xi(x) > u_i^* - b_i$ for all $i \in R$. Definition (2.18) implies that $k_i(x) = b_i$ for all $i \in R$. Moreover, in this case there exists at least one $i \in \{1,...,n\}$ for which $x_i > x_i^*$. Since $f_i$ is increasing on $[0, \mu_i]$ for $i = 1,...,n$ (a consequence of (H) and the fact that $\mu_i \leq \tilde{\delta}_i$), we conclude that there exists at least one $i \in \{1,...,n\}$ for which $f_i(x_i) > f_i(x_i^*)$. Consequently, we get from (3.12) and the fact that $u_i = k_i(x) = b_i$ for all $i \in R$:

$$\sum_{i=1}^{n}(n+1-i)u_i = \sum_{i \in R}(n+1-i)b_i + \sum_{i \notin R}(n+1-i)u_i^* \leq \min_{i=1,...,n}\left(((n-i)p_i + 1)f_i(x_i^*)\right) \leq \sum_{i=1}^{n}((n-i)p_i + 1)f_i(x_i).$$

Combining (3.17), (A.9) with the above inequality and using the fact $A \geq (1-L)^{-1} \sum_{i \in R} \sigma^i \gamma_i$, we conclude that (3.18) holds when $\Xi(x) > \tau$.

Case 2: $x \in S \setminus \Omega$, $d \in [0,1]^{n-1}$.

In this case, there exists at least one $i \in \{1,...,n\}$ for which $x_i > \mu_i$. Therefore, definition (2.19) implies $\Xi(x) > h = \min_{i=1,...,n}\left(\sigma^i(\mu_i - x_i^*)\right)$, and consequently definition (2.21) gives $P(x) = Q - \theta h$. Moreover, definition (2.18) gives $k_i(x) = b_i$ for all $i \in R$ (a direct consequence of the facts that $\tau < \tau^* \leq h$ and $\gamma_i = \tau^{-1}(u_i^* - b_i) \geq h^{-1}(u_i^* - b_i)$). Combining, we get from definition (2.20) and (3.17):

$$V(x^+) = \sum_{i=1}^{n} \sigma^i |x_i^+ - x_i^*| + A\Xi(x^+) + K \max\left(0, \sum_{i=1}^{n} I_i(x^+) - P(x^+)\right)$$
$$\leq \sum_{i=1}^{n} \sigma^i |x_i^+ - x_i^*| + K \max\left(0, \sum_{i=1}^{n} I_i(x^+) - Q + \theta h\right) + A\Xi(x^+) \quad (A.13)$$

Using (3.5), the facts that $u_i = k_i(x) = b_i$ for all $i \in R$, $Q - \theta h = \varepsilon \min_{i=1,...,n}\left((n+1-i)\mu_i\right)$ and $\sum_{i \in R}(n+1-i)b_i + \sum_{i \notin R}(n+1-i)u_i^* \leq \varepsilon C \min_{i=1,...,n}\left((n+1-i)\mu_i\right)$ (which both imply that $\sum_{i \in R}(n+1-i)b_i + \sum_{i \notin R}(n+1-i)u_i^* \leq C(Q - \theta h)$), we get:

$$\max\left(0, \sum_{i=1}^{n} I_i(x^+) - Q + \theta h\right) \leq (1-C)\max\left(0, \sum_{i=1}^{n} I_i(x) - Q + \theta h\right). \quad (A.14)$$

Combining (A.13) and (A.14), we get:

$$V(x^+) \le \sum_{i=1}^{n} \sigma^i \left|x_i^+ - x_i^*\right| + K(1-C)\max\left(0, \sum_{i=1}^{n} I_i(x) - Q + \theta h\right) + A\Xi(x^+). \quad (A.15)$$

Definition (2.19) in conjunction with (A.15) implies that the following inequality holds:

$$V(x^+) \le \sum_{i=1}^{n} \sigma^i \max\left(a_i - x_i^*, x_i^*\right) + K(1-C)\max\left(0, \sum_{i=1}^{n} I_i(x) - Q + \theta h\right) + A\sum_{i=1}^{n} \sigma^i \left(a_i - x_i^*\right). \quad (A.16)$$

The fact that there exists at least one $i \in \{1,...,n\}$ for which $x_i > \mu_i$, implies that

$$\sum_{i=1}^{n} I_i(x) = \sum_{i=1}^{n}(n+1-i)x_i \ge \min_{i=1,...,n}\left((n+1-i)\mu_i\right). \quad (A.17)$$

Using (A.16), (A.17) and the fact that $Q - \theta h = \varepsilon \min_{i=1,...,n}\left((n+1-i)\mu_i\right)$, we obtain:

$$V(x^+) \le \sum_{i=1}^{n} \sigma^i \max\left(a_i - x_i^*, x_i^*\right) + K\max\left(0, \sum_{i=1}^{n} I_i(x) - Q + \theta h\right)$$
$$+ A\sum_{i=1}^{n} \sigma^i\left(a_i - x_i^*\right) - KC(1-\varepsilon)\min_{i=1,...,n}\left((n+1-i)\mu_i\right) \quad (A.18)$$

Since $K \ge \dfrac{\sum_{i=1}^{n} \sigma^i \max\left(a_i - x_i^*, x_i^*\right) + A\sum_{i=1}^{n} \sigma^i\left(a_i - x_i^*\right) - (A+L)h}{(1-\varepsilon)C \min_{i=1,...,n}\left((n+1-i)\mu_i\right)}$, $\sum_{i=1}^{n} \sigma^i \left|x_i - x_i^*\right| \ge \Xi(x) > h$, we conclude from (A.18) and definition (2.20) that (3.18) holds. The proof is complete. ◁

**Proof of the Claim 2 made in the proof of Theorem 2.1:** Since $\sigma \in (0,1]$, we get for all $x \in S$:

$$\sigma^n \left|x - x^*\right| \le \sum_{i=1}^{n} \sigma^i \left|x_i - x_i^*\right| \le \left|x - x^*\right| \sum_{i=1}^{n} \sigma^i. \quad (A.19)$$

Similarly, using definition (2.19), we get for all $x \in S$:

$$0 \le \Xi(x) \le \sum_{i=1}^{n} \sigma^i \left|x_i - x_i^*\right| \le \left|x - x^*\right| \sum_{i=1}^{n} \sigma^i. \quad (A.20)$$

Using (A.20), the fact that $I_j(x) := \sum_{i=1}^{j} x_i$ for $j = 1,...,n$, definition (2.21) and the fact that $Q \ge \sum_{i=1}^{n} I_i(x^*)$ (a consequence of (3.5) and the fact that $Q \ge (1-C)\sum_{i=1}^{n} I_i(x^*) + \sum_{i=1}^{n}(n+1-i)u_i^*$), we get for all $x \in S$:

$$\max\left(0, \sum_{i=1}^{n} I_i(x) - P(x)\right) \leq \max\left(0, \sum_{i=1}^{n} I_i(x) - \sum_{i=1}^{n} I_i(x^*)\right) + \max\left(0, \sum_{i=1}^{n} I_i(x^*) - P(x)\right)$$

$$\leq \max\left(0, \sum_{i=1}^{n} (n+1-i)(x_i - x_i^*)\right) + \max\left(0, \sum_{i=1}^{n} I_i(x^*) - Q + \theta \min(h, \Xi(x))\right)$$

$$\leq \sum_{i=1}^{n} (n+1-i)|x_i - x_i^*| + \max\left(0, \sum_{i=1}^{n} I_i(x^*) - Q\right) + \theta \min(h, \Xi(x))$$

$$\leq \sum_{i=1}^{n} (n+1-i)|x_i - x_i^*| + \theta \Xi(x) \leq |x - x^*| \sum_{i=1}^{n} (n+1-i) + \theta |x - x^*| \sum_{i=1}^{n} \sigma^i$$

(A.21)

It follows from definition (2.20) and (A.19), (A.20), (A.21) that there exist constants $K_2 \geq K_1 > 0$ such that inequality (3.19) holds. The proof is complete. ◁